\documentclass{amsart}
\usepackage{amsmath, amsthm, amssymb}
\usepackage[all]{xy}

\theoremstyle{plain}
\newtheorem{theorem}{Theorem}[section]
\newtheorem{corollary}[theorem]{Corollary}
\newtheorem{lemma}[theorem]{Lemma}
\newtheorem{proposition}[theorem]{Proposition}

\theoremstyle{definition}

\newtheorem{remark}[theorem]{Remark}
\newtheorem{example}[theorem]{Example}

\theoremstyle{remark}

\begin{document}
\title[Homotopy theory of presheaves of simplicial groupoids]
      {Homotopy theory of presheaves\\
       of simplicial groupoids}
\author{Zhi-Ming Luo}
\address{Department of Mathematics\\
         University of Western Ontario\\
         London, Ontario, Canada N6A 5B7}
\email{zluo@uwo.ca}
\keywords{Homotopy category, presheaves of
simplicial groupoids, torsors}
\subjclass[2000]{Primary: 55U35;
Secondary: 55P15, 18F20}
\date{January 20, 2004}
\begin{abstract}
We show that the category of presheaves of simplicial groupoids on
a site $\mathcal{C}$ is a right proper simplicial model category.
We define $G$-torsor of presheaf of 2-groupoids $G$, presheaf of
simplicial groups $G$ and presheaf of simplicial groupoids $G$ on
a site $\mathcal{C}$ and classify $G$-torsors by the homotopy
classes $[\ast, \overline{W}G]$.
\end{abstract}

\maketitle

\section{Introduction}

The techniques of homotopy theory have been extremely fruitful in
other areas of mathematics such as algebraic K-theory \cite {J4}
and algebraic geometry \cite {M-V}. The goal of this paper is to
develop the machinery needed to do homotopy theory in category of
presheaves of simplicial groupoids on any small site
$\mathcal{C}$, and to connect this homotopy theory to the ordinary
homotopy theory of spaces.

Axiomatic homotopy theory is a natural extension of the ordinary
homotopy theory for topological spaces to other categories. It
comes from two systems of axioms. One is K.Brown's axioms for a
\emph{category of fibrant objects} \cite {BK}, \cite [Section
I.9.] {G-J}; the other one is Quillen's axioms for a \emph{closed
model category} \cite {Q1}, \cite {Q2}, \cite [Section II.1.]
{G-J}.

Quillen's axioms imply Brown's axioms. But on the level of sheaves
(presheaves), there are two homotopy theories. One is the local
theory which was constructed by Brown \cite {BK} for the category
of simplicial sheaves on a topological space, and by Jardine \cite
{J7}, \cite {J1} for the category of simplicial sheaves and
simplicial presheaves on any site. The other one is the global
theory for the corresponding categories, which was developed by
Brown and Gersten \cite {B-G}, Joyal \cite {Jo} and Jardine \cite
{J1}, respectively. The local theory satisfies Brown's axioms, the
global theory satisfies Quillen's axioms. The two theories are
distinct, since it is not true that every local fibration is a
global fibration \cite {J1}.

This paper is based on the Quillen closed model category axioms.

The central theorem of simplicial homotopy theory asserts that the
category {\bf S} of simplicial sets, equipped with three classes
of morphisms, namely cofibrations, fibrations and weak
equivalences, has a closed model structure \cite{Q1}.
Mathematicians have found many categories with closed model
structures. For example, the category of simplicial groupoids
(Dwyer-Kan \cite{D-K}, \cite [Section V.7.] {G-J}), the category
of 2-groupoids (Moerdijk-Svensson \cite{M-S}), the category of
simplicial presheaves (Jardine \cite{J1}), the category of
simplicial sheaves (Joyal \cite{Jo}) and so on. Crans \cite{Cr}
uses adjoint functors to prove that the categories of sheaves of
$2$-groupoids, of bisimplicial sheaves and of simplicial sheaves
of groupoids have closed model structures according to the
well-known closed model category of simplicial sheaves.

In Section 2 we use techniques based on Jardine's paper \cite{J1},
to prove that the categories of presheaves of simplicial groupoids
and of presheaves of 2-groupoids on any site $\mathcal{C}$ are
right proper simplicial model categories (Theorems \ref{T:main},
\ref{T:simpresg}, \ref{T:simpresgpro}, \ref{T:main2},
\ref{T:simpre2g} and \ref{T:2gpdpresgpro}). Both Crans \cite {Cr}
and Joyal-Tierney \cite {J-T5} obtain Quillen closed model
structures on the category of sheaves of simplicial groupoids in
the general sense (i.e., the simplicial groupoids are groupoid
objects in simplicial sets, not just groupoids enriched over
simplicial sets; in our case we restrict the simplicial groupoids
to be groupoids enriched in simplicial sets). They use the
classifying space functor $B$ to define the weak equivalences
whereas we use the functor $\overline{W}$.

Section 3 studies the classification of $G$-torsors for various
(pre)sheaves $G$. The definition of $G$-torsor and the relations
between $G$-torsors and the homotopy classes $[\ast, BG]$ are
well-known for sheaves of groups and groupoids $G$ (Sections 3.1
and 3.2). We extend the definition of $G$-torsors to presheaves of
2-groupoids, simplicial groups and simplicial groupoids $G$ on any
site $\mathcal{C}$ in Sections 3.3, 3.4 and 3.5. After we define
the category $Tor(\ast, G)$ of all $G$-torsors we obtain similar
classifying results for $G$-torsors: $\pi_0Tor(\ast, G) \cong
[\ast, \overline{W}G]$ (Theorem \ref{T:2gtor}, Corollary
\ref{C:Gsheaf} and Theorem \ref{T:maintor}). All these definitions
and results for (pre)sheaves of groups, groupoids, 2-groupoids and
simplicial groups are special cases of those for presheaves of
simplicial groupoids. Joyal-Tierney \cite {J-T6} obtain similar
results for sheaves of simplicial groupoids $G$, but their
definition of $G$-torsor is different from ours; our definition is
much more flexible.

\section{Presheaves of simplicial groupoids}

\subsection{Presheaves of simplicial groupoids}

Suppose that $C$ is a simplicial object in the category of small
categories. Write $E_C$ for the category constructed using the
following variant of the Grothendieck construction.

The set of objects of $E_C$ consists of all pairs $(x,n)$ with
$x\in C_n$. A morphism $(f,\theta): (x,m) \to (y,n)$ is a pair
consisting of an ordinal number map $\theta: {\bf m} \to {\bf n}$
and a morphism $f: x \to \theta^\ast y$ of $C_m$. Composition is
defined in the natural way. There is an obvious forgetful functor
$\pi: E_C \to \Delta$ which takes values in the ordinal number
category $\Delta$.

The segment category $Seg({\bf n})$ of subintervals $[j,n]$ of
${\bf n} = [0,n]$ has as objects the intervals $[j, n] = \{ j,
j+1, ..., n \}$ and as morphisms the inclusions of intervals.
$Seg({\bf n})$ can be identified with the opposite ${\bf n}^{op}$
via the functor $[j,n] \mapsto j$. There is a functor $c_n: {\bf
n}^{op} \to \Delta$ which is defined on objects by $j \mapsto {\bf
n-j}$, and which sends morphisms to inclusions.

An $n$-cocycle taking values in the simplicial category $C$ is a
functor $X: {\bf n}^{op} \to E_C$ which is a lifting of $c_n$ in
the sense that the diagram of functors
$$
\xymatrix{&E_C \ar[d]^-\pi\\
{\bf n}^{op}\ar[ur]^-X \ar[r]_-{c_n}&\Delta}
$$
commutes. This is a generalization of the definition of an
$n$-cocycle taking values in a groupoid enriched in simplicial
sets, in view of the identification of the categories $Seg({\bf
n})$ and ${\bf n}^{op}$.

The $n$-cocycle $X: {\bf n}^{op} \to E_C$ is otherwise described
as a string of arrows
$$(x_0,n)\leftarrow (x_1, n-1) \leftarrow \cdots \leftarrow (x_n,
0)$$ each of which has the form $(\alpha_i, d^0)$, with $\alpha_i:
x_{n-i} \to d_0(x_{n-i-1})$. This means that the string consists
of objects $x_i \in C_{n-i}$ and morphisms $x_i \to d_0(x_{i-1})$.

Every ordinal number map $\theta: {\bf m} \to {\bf n}$ induces a
commutative diagram
$$
\xymatrix{{\bf
m-i}\ar[d]_-{\theta_i}\ar[r]^-\cong&[i,m]\ar[d]_-{\theta_i}\ar[r]&{\bf
m}\ar[d]^-\theta\\
{\bf n}- \theta({\bf i}) \ar[r]_-\cong&[\theta(i),n] \ar[r] &{\bf
n}}
$$
and there is a corresponding diagram
$$
\xymatrix{(\theta^\ast_0 x_{\theta(0)},m) \ar[d]_-{(1,\theta_0)}&
(\theta^\ast_1 x_{\theta(1)},m-1) \ar[d]_-{(1,\theta_1)}\ar[l] &
\cdots\ar[l] & (\theta^\ast_m x_{\theta(m)},0)
\ar[d]_-{(1,\theta_m)}\ar[l]\\
(x_{\theta(0)}, n-\theta(0))& (x_{\theta(1)},
n-\theta(1))\ar[l]&\cdots\ar[l]&(x_{\theta(m)},
n-\theta(m))\ar[l]}
$$
The string on top is denoted by $\theta^\ast X$.

In this way, a simplicial set $ \overline{W}C$ is defined, with $
\overline{W}C_n$ given by the set of $n$-cocycles in $C$. The
functoriality follows from the relations
$$\theta_{\tau(i)}\tau(i) = (\theta\tau)_i$$
associated to composeable ordinal number maps
$$\xymatrix@1{{\bf k}\ar[r]^-\tau&{\bf m}\ar[r]^-\theta&{\bf n}}.$$

There is an obvious function
$$j: dBC_n = (BC_n)_n \to \overline{W}C_n$$
which sends a string
$$\xymatrix@1{x_0 &x_1 \ar[l]_-{\alpha_n} & \cdots
\ar[l]_-{\alpha_{n-1}} & x_n \ar[l]_-{\alpha_1}}$$ to the cocycle
consisting of the object $x_0 \in C_n$ and the morphisms
$$d^j_0\alpha_{n-j}: d^j_0x_j \to d^j_0x_{j-1},$$
or rather to the string
$$ (x_0, n) \leftarrow (d_0x_1, n-1) \leftarrow \cdots \leftarrow
(d^n_0x_n, 0)$$ in the Grothendieck construction $E_C$.

Suppose that $\theta: {\bf m} \to {\bf n}$ is an ordinal number
map. One checks that the composite
$$\xymatrix@1{dBC_n \ar[r]^-j & \overline{W}C_n
\ar[r]^-{\theta^\ast} & \overline{W}C_m}$$ sends the string of
arrows $x_0 \leftarrow x_1 \leftarrow \cdots \leftarrow x_n$ in
$C_n$ to the string
$$(\theta^\ast_0d^{\theta(0)}_0x_{\theta(0)}, m) \leftarrow
(\theta^\ast_1d^{\theta(1)}_0x_{\theta(1)}, m-1) \leftarrow \cdots
\leftarrow (\theta^\ast_md^{\theta(m)}_0x_{\theta(m)}, 0)$$ while
the composite
$$\xymatrix@1{dBC_n \ar[r]^-{\theta^\ast} &dBC_m \ar[r]^-j &
\overline{W}C_m}$$ sends that same string in $C_n$ to the string
$$(\theta^\ast x_{\theta(0)}, m) \leftarrow (d_0\theta^\ast x_{\theta(1)}, m-1)
\leftarrow \cdots \leftarrow (d^m_0\theta^\ast x_{\theta(m)},
0).$$ Since $\theta^\ast_i d^{\theta(i)}_0 x_{\theta(i)} =
d^i_0\theta^\ast x_{\theta(i)}$ (from above diagram of $[i,m],
[\theta(i),n]$ and {\bf m}, {\bf n}), it follows that the maps $j$
respect the simplicial structures.

A simplicial groupoid $G$ (or simplicial group $G$) is a
simplicial object in the category of small groupoids (or groups),
so the functor $ \overline{W}$ can act on it and produce a
simplicial set $ \overline{W}G$ (see \cite [Section V.4 and V.7]
{G-J}).

\begin{lemma}\label{L:dBW}
Suppose that $G$ is a simplicial groupoid (groupoid enriched in
simplicial sets). Then the map $j: dBG \to \overline{W}G$ is a
weak equivalence.
\end{lemma}
\begin{proof}
Since both functors $dB$ and $ \overline{W}$ are right adjoint
functors, they preserve limits and products, so they preserve
homotopy equivalences as well \cite[pp.303,304]{G-J}. They also
preserve disjoint unions. If $H$ is a simplicial group, the map
$j: dBH \to \overline{W}H$ classifies the $H$-bundle $dEH \to
dBH$, and so $j$ is a weak equivalence for simplicial groups.
Every simplicial groupoid $G$ is homotopy equivalent to a disjoint
union of simplicial groups.
\end{proof}

\begin{lemma}
Suppose that $H$ is a simplicial groupoid. Then the following
statements hold:
\begin{itemize}
\item[1)] $dBH$ is a Kan complex.
\item[2)] $\overline{W}H$ is a Kan complex.
\end{itemize}
\end{lemma}
\begin{proof}
Every simplicial groupoid $H$ contains a strong deformation
retraction $\bigsqcup_x H(x,x)$, where the vertices $x$ are
indexed over a set of representatives of $\pi_0H$. Suppose that
$X$ denotes either $dB$ or $\overline{W}$. Then every lifting
problem
$$\xymatrix{
\wedge^n_k \ar[r]\ar[d] &X(H)\\
\triangle^n\ar@{.>}[ur]}
$$
is isomorphic to a lifting problem
$$\xymatrix{
\wedge^n_k \ar[r]\ar[d] &\bigsqcup_x X(H(x,x))\\
\triangle^n\ar@{.>}[ur]}
$$
It follows that $X(H)$ is a Kan complex for all simplicial
groupoids $H$ if and only if $X(G)$ is a Kan complex for all
simplicial groups $G$.

In the first case, the bisimplicial set $BG$ consists of connected
simplicial sets $BG_n$ in each horizontal degree $n$, and
therefore satisfies the $\pi_\ast$-Kan condition \cite [Lemma
IV.4.2.] {G-J}. In vertical degree $k$, $BG_{\ast, k} = G^{\times
k}$, which is a Kan complex since all simplicial groups are Kan
complexes. It follows from Lemma IV.4.8 of \cite {G-J} that $dBG$
is a Kan complex.

For the second statement, we know that $\overline{W}G$ is a Kan
complex if $G$ is a simplicial group, by standard theory \cite
[Corollary V.6.8.] {G-J}.
\end{proof}

\bigskip
Let $ \mathcal{C}$ be a fixed small Grothendieck site. $s${\bf
Gd}Pre($\mathcal{C}$) is the category of presheaves of simplicial
groupoids on $ \mathcal{C}$; its objects are the contravariant
functors from $ \mathcal{C}$ to the category $s${\bf Gd} of
simplicial groupoids, and its morphisms are natural
transformations.

If we see a presheaf of simplicial groupoids $G$ as a simplicial
object in the category of presheaves of groupoids, then $
\overline{W}G$ is a simplicial presheaf. That means there is a
functor $ \overline{W}: s{\bf Gd}Pre(\mathcal{C}) \to {\bf S
}Pre(\mathcal{C})$.

Recall the adjunction between the loop groupoid functor $G: {\bf
S} \to s{\bf Gd} $ and the universal cocycle functor
$\overline{W}$ \cite[Lemma V.7.7]{G-J}. By applying these functors
sectionwise to simplicial presheaves and presheaves of simplicial
groupoids, one obtains functors
$$G: {\bf S}Pre(\mathcal{C})\rightleftarrows s{\bf Gd}Pre(\mathcal{C}): \overline{W}$$

So there is

\begin{proposition}\label{P:adjoint}
The functor $G: {\bf S}Pre(\mathcal{C}) \to s{\bf Gd}Pre(
\mathcal{C})$ is left adjoint to the functor $\overline{W}$.
\end{proposition}

\bigskip
We recall the definition of a closed model category structure on a
category $\mathcal{D}$. Closed model categories are an abstract
setting in which to do homotopy theory.

A \emph{Quillen closed model category $\mathcal{D}$} is a category
which is equipped with three classes of morphisms, called
cofibrations, fibrations and weak equivalences which together
satisfy the following axioms \cite{Q1}, \cite{Q2}, \cite{G-J}:
  \begin{itemize}
     \item[{\bf CM1:}] The category $\mathcal{D}$ is closed under all finite
            limits and colimits.
     \item[{\bf CM2:}] Suppose that the following diagram
            commutes in $\mathcal{D}$:
         $$\xymatrix{ X \ar[rr]^-g \ar[dr]_-h && Y\ar[dl]^-f\\
                        &   Z }
          $$
        If any two of $f,g$ and $h$ are weak equivalences, then so is the
        third.
     \item[{\bf CM3:}] If $f$ is a retract of $g$ and $g$ is a weak
         equivalence, fibration or cofibration, then so is $f$.
      \item[{\bf CM4:}] Suppose that we are given a commutative diagram
         $$\xymatrix{ U\ar[r]\ar[d]_i   & X \ar[d]^p\\
             V \ar[r]\ar@{.>}[ur]    & Y  }
            $$
        where $i$ is a cofibration and $p$ is a fibration. Then the
        lifting exists, making the diagram commute, if either $i$ or $p$
        is also a weak equivalence.
     \item[{\bf CM5:}] Any map $f: X \to Y $ may be factored:
        \begin{itemize}
           \item[(a)] $f = p \cdot i $ where $p$ is a fibration and $i$ is a trivial
                  cofibration, and
            \item[(b)] $f = q \cdot j $ where $q$ is a trivial fibration and $j$ is a
              cofibration.
         \end{itemize}
       \end{itemize}

An object $X$ is called \emph{cofibrant} if the map from the
initial object $\emptyset$, to $X$ is a cofibration. An object $X$
is called \emph{fibrant} if the map from $X$ to the final object
$\ast$, is a fibration. The category obtained from $\mathcal{D}$
by formally inverting the weak equivalences is called the
\emph{homotopy category} associated to $\mathcal{D}$, and denoted
$Ho(\mathcal{D})$.

A category $ \mathcal{D}$ is a \emph{simplicial category} if there
is a mapping space functor
$${\bf Hom}_{ \mathcal{D}}(\cdot, \cdot):
\mathcal{D}^{op} \times \mathcal{D } \to \bf{S}$$ with the
properties that for $A$ and $B$ objects in $ \mathcal{D}$
  \begin{itemize}
   \item[(1)] ${\bf Hom}_\mathcal{D}(A,B)_0 = \mbox{hom}_
   \mathcal{D}(A,B)$;
   \item[(2)] the functor ${\bf Hom}_{ \mathcal{D}}(A, \cdot):
    \mathcal{D } \to \bf{S}$ has a left adjoint
    $$A \otimes \cdot: {\bf S} \to \mathcal{D}$$
    which is associative in the sense that there is an isomorphism
    $$A\otimes(K \times L) \cong (A\otimes K) \otimes L,$$
    natural in $A \in \mathcal{D}$ and $K,L \in {\bf S}$;
    \item[(3)] The functor ${\bf Hom}_{ \mathcal{D}}(\cdot, B):
    \mathcal{D}^{op}  \to \bf{S}$ has a left adjoint
    $${\bf hom}_{ \mathcal{D}}(\cdot, B): {\bf S} \to
    \mathcal{D}^{op}.  $$
   \end{itemize}

A \emph{simplicial model category} $ \mathcal{D}$ is both a closed
model category and a simplicial category which satisfies the
following axiom:
  \begin{itemize}
  \item[{\bf SM7}] Suppose $j: A \to B$ is a cofibration and $q: X
  \to Y$ is a fibration. Then
  $$\xymatrix@1{{\bf Hom}_\mathcal{D}(B,X)\ar[r]^-{(j^\ast,q_\ast)}
  &{\bf Hom}_\mathcal{D}(A,X) \times_{{\bf Hom}_\mathcal{D}(A,Y)} {\bf
  Hom}_\mathcal{D}(B,Y)}$$
  is a fibration of simplicial sets, which is trivial if $j$ or
  $q$ is trivial.
 \end{itemize}

A \emph{right proper closed model category $\mathcal{D}$} is a
closed model category such that: {\bf P1} the class of weak
equivalences is closed under base change by fibrations. In other
words, axiom {\bf P1} says that, given a pullback diagram
$$\xymatrix{
X \ar[r]^-{g_\ast} \ar[d] & Y\ar[d]^-p\\
Z \ar[r]_-g& W}
$$
of $\mathcal{D}$ with $p$ a fibration, if $g$ is a weak
equivalence then so is $g_\ast$.

\bigskip
The category {\bf S}Pre($ \mathcal{C}$) of simplicial presheaves
has a Quillen closed model structure \cite {J1}, hence we can give
some definitions in the category $s${\bf Gd}Pre($ \mathcal{C}$).

A map $f: X \to Y$ in the category $s${\bf Gd}Pre($ \mathcal{C}$)
is said to be a \emph{fibration} if the induced map
$\overline{W}(f): \overline{W}X \to \overline{W}Y$ is a global
fibration in the category {\bf S}Pre($ \mathcal{C}$) in the sense
of \cite{J1}.

A map $g: Z \to U$ in the category $s${\bf Gd}Pre($ \mathcal{C}$)
is said to be a \emph{weak equivalence} if the induced map
$\overline{W}(g): \overline{W}Z \to \overline{W}U$ is a
topological weak equivalence in the category {\bf S}Pre($
\mathcal{C}$) in the sense of \cite{J1}.

A \emph{cofibration} in the category $s${\bf Gd}Pre($
\mathcal{C}$) is a map of presheaves of simplicial groupoids which
has the left lifting property with respect to all both fibrations
and weak equivalences.

\begin{theorem}\label{T:main}
The category $s${\bf Gd}Pre($\mathcal{C}$), with the classes of
fibrations, weak equivalences and cofibrations as defined above,
satisfies the axioms for a closed model category.
\end{theorem}
\begin{proof}
See \cite {Luo}.
\end{proof}

\begin{remark}
Crans \cite {Cr} and Joyal-Tierney \cite {J-T5} provide two
different Quillen closed model structures on the category of
sheaves of simplicial groupoids in the general sense (i.e., the
simplicial groupoids are groupoid objects in simplicial sets, not
just groupoids enriched) (see the comment in the introduction of
\cite {J-T5}). They use the classifying space functor $B$ to
define the weak equivalences whereas we use the functor
$\overline{W}$.
\end{remark}

Let $X$ be a presheaf of simplicial groupoids, and let $K$ be a
simplicial set. The presheaf of simplicial groupoids $X \otimes K$
is defined at $U \in \mathcal{C}$ by
$$X \otimes K (U) = X(U) \otimes K, $$
The presheaf of simplicial groupoids $X^K$ is defined at $U \in
\mathcal{C}$ by
$$X^K(U) = {\bf hom}(K, X(U)),$$
where ${\bf hom}(K, X(U))$ is a simplicial groupoid.

For presheaves of simplicial groupoids $G$ and $H$, define a
simplicial set ${\bf Hom}(G,H)$ by requiring that the
$n$-simplices be maps of presheaves of simplicial groupoids of the
form $G \otimes \triangle^n \to H$.

Let $i: K \to L$ be a cofibration in {\bf S} and $q: U \to V$ be a
fibration in $s${\bf Gd}Pre($ \mathcal{C}$), then the map
$$U^L \stackrel{(q_*,i^*)}{\longrightarrow} V^L \times_{V^K} U^K$$
is a fibration, which is trivial if either $i$ or $q$ is trivial,
since when we apply the functor $ \overline{W}$ we get a similar
map in {\bf S}Pre($ \mathcal{C}$) and the category {\bf S}Pre($
\mathcal{C}$) is a simplicial model category.

Suppose that $p: G \to H$ is a fibration and that $i: X \to Y$ is
a cofibration of presheaves of simplicial groupoids, then
$$ {\bf Hom}(Y, G) \stackrel{(i^*, p_*)}{\longrightarrow}
{\bf Hom}(X, G) \times_{{\bf Hom}(X, H)} {\bf Hom}(Y, H)$$ is a
fibration of simplicial sets, which is trivial if $i$ or $p$ is
trivial by \cite [Proposition II.3.13] {G-J}.

The mapping space functor satisfies the axiom {\bf SM7}, so there
is

\begin{theorem}\label{T:simpresg}
The category $s${\bf Gd}Pre($ \mathcal{C}$) is a simplicial model
category.
\end{theorem}

\begin{theorem}\label{T:simpresgpro}
The category $s${\bf Gd}Pre($ \mathcal{C}$) is right proper.
\end{theorem}
\begin{proof}
Given a pullback diagram in $s${\bf Gd}Pre($ \mathcal{C}$)
$$\xymatrix{
X \ar[r]^-{g_\ast} \ar[d] & Y\ar[d]^-p\\
Z \ar[r]_-g& W}
$$
with $p$ a fibration and $g$ a weak equivalence, there exists a
pullback diagram in {\bf S}Pre($ \mathcal{C}$)
$$\xymatrix{
\overline{W}X \ar[r]^-{\overline{W}g_\ast} \ar[d] & \overline{W}Y\ar[d]^-{\overline{W}p}\\
\overline{W}Z \ar[r]_-{\overline{W}g}& \overline{W}W}
$$
$\overline{W}$ preserves fibrations and weak equivalences, hence
$\overline{W}p$ is a fibration and $\overline{W}g$ is a weak
equivalence in {\bf S}Pre($ \mathcal{C}$). {\bf S}Pre($
\mathcal{C}$) is proper, so $\overline{W}g_\ast$ is a weak
equivalence as well, hence the map $g_\ast$ is a weak equivalence.
So the axiom {\bf P1} holds.
\end{proof}

\subsection{Presheaves of 2-groupoids}

{\bf 2-Gpd}Pre($\mathcal{C}$) is the category of presheaves of
2-groupoids on $\mathcal{C}$; its objects are the contravariant
functors from $\mathcal{C}$ to the category {\bf 2-Gpd} of
2-groupoids, and its morphisms are natural transformations.

Recall the adjunction \cite {Luo1}:
$$\pi: s{\bf Gd} \rightleftarrows {\bf 2-Gpd}: B$$
By applying these functors sectionwise to presheaves of simplicial
groupoids and presheaves of 2-groupoids, one obtains functors
$$\pi: s{\bf Gd}Pre(\mathcal{C}) \rightleftarrows {\bf 2-Gpd}Pre(\mathcal{C}): B$$
and there is

\begin{proposition}\label{P:adj}
The functor $\pi: s{\bf Gd}Pre(\mathcal{C}) \to {\bf
2-Gpd}Pre(\mathcal{C})$ is left adjoint to the functor $B$.
\end{proposition}

The category $s${\bf Gd}Pre($ \mathcal{C}$) of presheaves of
simplicial groupoids has a Quillen closed model structure (see
Section 2.1), hence we can give some definitions in the category
{\bf 2-Gpd}Pre($ \mathcal{C}$).

A map $f: X \to Y$ in the category {\bf 2-Gpd}Pre($ \mathcal{C}$)
is said to be a \emph{fibration} if the induced map $B(f): BX \to
BY$ is a fibration in the category $s${\bf Gd}Pre($ \mathcal{C}$).

A map $g: Z \to U$ in the category {\bf 2-Gpd}Pre($ \mathcal{C}$)
is said to be a \emph{weak equivalence} if the induced map $B(g):
BZ \to BU$ is a weak equivalence in the category $s${\bf Gd}Pre($
\mathcal{C}$).

A \emph{cofibration} in the category {\bf 2-Gpd}Pre($
\mathcal{C}$) is a map of presheaves of 2-groupoids which has the
left lifting property with respect to all fibrations and weak
equivalences.

\begin{theorem}\label{T:main2}
The category {\bf 2-Gpd}Pre($\mathcal{C})$, with the classes of
fibrations, weak equivalences and cofibrations as defined above,
satisfies the axioms for a closed model category.
\end{theorem}
\begin{proof}
See \cite {Luo1} or \cite {Luo}.
\end{proof}

Let $X$ be a presheaf of 2-groupoids, and let $K$ be a simplicial
set. The presheaf of 2-groupoids $X \otimes K$ is defined at $U
\in \mathcal{C}$ by
$$X \otimes K (U) = X(U) \otimes K, $$
The presheaf of 2-groupoids $X^K$ is defined at $U \in
\mathcal{C}$ by
$$X^K(U) = (X(U))^K.$$

For presheaves of 2-groupoids $G$ and $H$, define a simplicial set
${\bf Hom}(G,H)$ by requiring that the $n$-simplices be maps of
presheaves of 2-groupoids of the form $G \otimes \triangle^n \to
H$.

Let $i: K \to L$ be a cofibration in {\bf S} and $q: U \to V$ be a
fibration in {\bf 2-Gpd}Pre($ \mathcal{C}$), then the map
$$U^L \stackrel{(q_*,i^*)}{\longrightarrow} V^L \times_{V^K} U^K$$
is a fibration, which is trivial if either $i$ or $q$ is trivial,
since apply the functor $B$ to get a similar map in $s${\bf
Gd}Pre($ \mathcal{C}$) and the category $s${\bf Gd}Pre($
\mathcal{C}$) is a simplicial model category.

Suppose that $p: G \to H$ is a fibration and that $i: X \to Y$ is
a cofibration of presheaves of 2-groupoids, then
$$ {\bf Hom}(Y, G) \stackrel{(i^*, p_*)}{\longrightarrow}
{\bf Hom}(X, G) \times_{{\bf Hom}(X, H)} {\bf Hom}(Y, H)$$ is a
fibration of simplicial sets, which is trivial if $i$ or $p$ is
trivial by \cite [Proposition II.3.13] {G-J}.

The mapping space functor satisfies the axiom {\bf SM7}, so there
is

\begin{theorem}\label{T:simpre2g}
The category {\bf 2-Gpd}Pre($ \mathcal{C}$) is a simplicial model
category.
\end{theorem}

\begin{theorem}\label{T:2gpdpresgpro}
The category {\bf 2-Gpd}Pre($ \mathcal{C}$) is right proper.
\end{theorem}
\begin{proof}
The functor $B$ preserves pullbacks and the category $s${\bf
Gd}Pre($\mathcal{C}$) is right proper
(Theorem~\ref{T:simpresgpro}).
\end{proof}

\section{Classification of torsors}


\subsection{Torsors for sheaves of groups }

Suppose that $G$ is a sheaf of groups on $\mathcal{C}$. A (right)
\emph{G-torsor} is a non-empty sheaf $E$ (meaning $E \to 1 $ is
surjective) equipped with a free (right) $G$-action $a: E \times G
\to E$ which is transitive \cite {J-T}, \cite {J5}, \cite {M-M}.

The non-abelian cohomology object $H^1 (\mathcal{C}, G)$ is the
set of isomorphism classes of $G$-torsors on the site
$\mathcal{C}$, as usual.

It is shown in \cite {J3}, \cite {J5} :

\begin{theorem}\label{T:gptor}
There is an isomorphism
$$[\ast, BG] \cong H^1 (\mathcal{C}, G)$$
for any sheaf of groups $G$ on any Grothendieck site
$\mathcal{C}$.
\end{theorem}

\begin{remark}
The sheaf of groupoids Tor($G$) consisting of groupoids
Tor$(G|_U)$ of all $G|_U$-torsors over $U \in \mathcal{C}$ is a
canonical gerbe \cite {G}, \cite {B}. If $ \mathcal{C}$ has a
terminal object $\ast$ and $ \mathcal{G}$ is a $G$-gerbe such that
$\mathcal{G}(\ast) \neq \emptyset$, then the $G$-torsors are
bijective to the objects of the groupoid $\mathcal{G}(\ast)$ up to
a unique isomorphism \cite [Proposition 5.2.5] {Bry}, i.e., the
isomorphism classes of $G$-torsors are bijective to the
isomorphism classes of objects of groupoid $ \mathcal{G}(\ast)$,
so $[\ast, BG]$ is also in a natural bijection with the set of
isomorphism classes of objects of $ \mathcal{G}(\ast)$.
\end{remark}

\begin{remark}
The central theorem of non-abelian cohomology asserts that the
second cohomology group $H^2( \mathcal{C}, A)$ of the sheaf of
abelian groups $A$ on the site $ \mathcal{C}$ is isomorphic to the
group of equivalence classes of $A$-gerbes over $ \mathcal{C}$
\cite {G}, \cite {Bry}.
\end{remark}


\subsection{Torsors for sheaves of groupoids}

Let {\bf $\mathcal{E}$} denote the topos $Shv(\mathcal{C})$ of
sheaves on the site $\mathcal{C}$. We can see a sheaf of groupoids
$G$ as a reflexive graph in {\bf $\mathcal{E}$} \cite {J-T}:
$$\xymatrix{G_0 \ar[r]^u & G_1 \ar@<+0.5ex>[r]^s \ar@<-0.5ex>[r]_t & G_0},$$
where $su = tu = id $, provided with an associative composition
$c: G_1 \times_{G_0} G_1 \to G_1$, for which the elements of $G_0$
are units (via $u$), and each element of $G_1$ is invertible.

There are two kinds of definitions of $G$-torsor of sheaves of
groupoids $G$. One is in the sense of \cite {J-T}, the other one
is in the sense of \cite {J5}. We can show that they are
equivalent. The definitions of $G$-torsors in following sections
are based on the definition in \cite {J5}.

A (right) \emph{G-torsor} in the sense of \cite {J-T} is a
non-empty sheaf $E$ over $G_0$ in {\bf $\mathcal{E}$} equipped
with a free (right) action $a: E \times_{G_0} G_1 \to E$ which is
transitive, where $E \times_{G_0} G_1$ is the pullback:
$$\xymatrix{
E \times_{G_0} G_1\ar[r]^-{\pi_1} \ar[d]_{\pi_2} & E\ar[d]^f\\
G_1\ar[r]_t
 & G_0}$$
In other words, if we denote by $\mathbb{E}$ the groupoid given by
the top of the diagram (cf.~\cite {J-T5}):
$$\xymatrix{
E \times _{G_0}G_1 \ar@<+0.5ex>[r]^-a \ar@<-0.5ex>[r]_-{\pi_1} \ar[d]_-{\pi_2}& E \ar[d]^-f\\
G_1 \ar@<+0.5ex>[r]^-s \ar@<-0.5ex>[r]_-t & G_0}$$ then the
$G$-torsor $E$ is a sheaf such that the sheaf of groupoids
$\mathbb{E}$ is trivial and locally connected.

A \emph{G-torsor} in the sense of \cite {J5} is a simplicial sheaf
map $Y \to BG$ such that $Y$ is a homotopy colimit of some functor
$X$ and the canonical map $Y \to \ast$ is a local weak
equivalence.

We can construct a sequence of pullbacks
$$\xymatrix{
\cdot\cdot\cdot\ar[r]&E \times_{G_0} G_1 \times_{G_0}
G_1\cdot\cdot\cdot\times_{G_0} G_1\ar[r] \ar[d]
&\cdot\cdot\cdot\ar[r]&E \times_{G_0} G_1 \times_{G_0} G_1\ar[r]
\ar[d]&E \times_{G_0} G_1\ar[r]^-{\pi_1} \ar[d]_{\pi_2}
& E\ar[d]^f\\
\cdot\cdot\cdot\ar[r]& G_1 \times_{G_0}
G_1\cdot\cdot\cdot\times_{G_0} G_1\ar[r] &\cdot\cdot\cdot\ar[r]&
G_1 \times_{G_0} G_1\ar[r]& G_1\ar[r]_t & G_0}$$ Let $Y$ be the
simplicial sheaf with sheaf $Y_n = E \times_{G_0} G_1 \times_{G_0}
G_1\cdot\cdot\cdot\times_{G_0} G_1$ ($n$ factors of $G_1$, so $Y_0
= E$) and the natural structure maps. Thus the above diagram is a
simplicial sheaf map $\pi: Y \to BG$ such that all diagrams
$$\xymatrix{
Y_n\ar[r]^{0^\ast} \ar[d]_{\pi} & Y_0\ar[d]^f\\
BG_n\ar[r]_{0^\ast}
 & BG_0}$$
are pullbacks, where $0^\ast$ is the map corresponding to the
ordinal number map $0: {\bf 0 \to n}$ which picks out the object
$0$.

So $Y$ has the homotopy colimit structure $Y\cong EX$ for some
sheaf-valued functor $X: G \to {\bf Shv(\mathcal{C})}$ defined on
the small sheaf of groupoids $G$ over the site $\mathcal{C}$ \cite
{J5}. Since $Y$ is obtained from $BG$ by pullbacks, $Y$ is the
nerve $EX$ of the sheaf of groupoids:
$$\xymatrix{E \ar[r]^-u & E \times_{G_0} G_1 \ar@<+0.5ex>[r]^-a \ar@<-0.5ex>[r]_-{\pi_1} & E}$$
where $u: e \to (e, id_{f(e)})$. Since $G$ acts transitively on
the sheaf $Y_0 = E$, $Y$ has only one component; $G$ acts freely
on the sheaf $Y_0 = E$, so there is only an identity map from
every object in $Y$ to itself. So the canonical map $Y \to \ast$
is a local weak equivalence.

Every $G$-torsor $E$ in the sense of \cite {J-T} uniquely
determines a simplicial sheaf map $Y \to BG$ such that the map
satisfies the above pullback condition and the canonical map $Y
\to \ast$ is a locally weak equivalence. The simplicial sheaf map
$Y \to BG$ is just a $G$-torsor in the sense of \cite {J5}.

On the other hand, for each $G$-torsor $Y \to BG$ in the sense of
\cite {J5} the simplicial sheaf $Y$ is formed by pullbacks
$$\xymatrix{
Y_n\ar[r]^{0^\ast} \ar[d]_{\pi} & Y_0\ar[d]^f\\
BG_n\ar[r]_{0^\ast}
 & BG_0}$$
such that the diagram
$$\xymatrix{
Y_0 \times _{G_0}G_1 \ar@<+0.5ex>[r]^-{d_1 = s} \ar@<-0.5ex>[r]_-{0^\ast = d_0 = t} \ar[d]_-{\pi_2}& Y_0 \ar[d]^-f\\
G_1 \ar@<+0.5ex>[r]^-{d_1 = s} \ar@<-0.5ex>[r]_-{0^\ast = d_0 = t}
& G_0}$$ commutes, hence $d_1: Y_0 \times_{G_0} G_1 \to Y_0$ is an
action. $Y \to \ast$ is a local weak equivalence, so the action is
free and transitive and $Y_0$ is a $G$-torsor in the sense of
\cite {J-T}.

It's obvious that there is a bijection between sets of the two
kinds of $G$-torsors. So the definitions of \emph{G-torsor} in the
sense of \cite {J-T} and \cite {J5} are equivalent.

Take the special case $G_0 = \ast$, then the sheaf of groups $G$
is a special sheaf of groupoids
$$\xymatrix{\ast \ar[r]^-u & G \ar@<+0.5ex>[r]^-s \ar@<-0.5ex>[r]_-t & \ast}$$
A $G$-torsor is a $G$-space $Y_0$ such that
$$Y \cong EG \times_G Y_0 \simeq \ast$$
which is equivalent to the above definition.

The set of $G$-torsors and natural transformations between them
forms a category which will be denoted by Tors($\ast, G$). This
category is a groupoid because every morphism of $G$-torsors is an
isomorphism \cite {J5}. Denote the path components of Tors($\ast,
G$) by $\pi_0$Tors($\ast, G$).

It is shown in \cite{J5}

\begin{theorem}\label{T:gpdtor}
The function $\pi_0Tors(\ast, G) \to [\ast, BG]$ is a natural
bijection for any sheaf of groupoids $G$ on any Grothendieck site
$\mathcal{C}$.
\end{theorem}

\begin{remark}
As pointed out by Jardine \cite {J6}, there is a mistake in the
proof of \cite {J5}. Jardine corrected it in \cite {J6}. The
correction will also appear in the proof of Theorem
\ref{T:maintor}.
\end{remark}

\begin{remark}
Let $H( \mathcal{C}, G)$ denote the category of $G$-torsors and
$H^1 ( \mathcal{C}, G)$ the set of connected components of $H(
\mathcal{C}, G)$ \cite {J-T}. Thus Theorem \ref{T:gpdtor} has same
form as Theorem \ref{T:gptor}.
\end{remark}


\subsection{Torsors for presheaves of 2-groupoids}

Let {\bf 2-Sets} denote the 2-category in which the objects are
the sets, the morphisms are the maps between sets, the 2-arrows
are the maps $X \times I \to Y$ where $X$ and $Y$ are two sets and
$I$ is the 2-point set $\{0, 1\}$. So there is a unique 2-arrow
between each pair of parallel 1-arrows. In the language of
internal categories, {\bf 2-Sets} is a category enriched in
categories. Let ${\bf 2-Sets}_0$, ${\bf 2-Sets}_1$ and ${\bf
2-Sets}_2$ denote the classes of objects, 1-arrows and 2-arrows in
{\bf 2-Sets}, respectively. The underlying category
$\xymatrix@1{{\bf 2-Sets}_1 \ar@<-1ex>[r] \ar@<+1ex>[r] & {\bf
2-Sets}_0 \ar[l]}$ of {\bf 2-Sets} is just the ordinary category
{\bf Sets} of sets, $\xymatrix@1{{\bf 2-Sets}_2 \ar@<-1ex>[r]
\ar@<+1ex>[r] & {\bf 2-Sets}_0 \ar[l]}$ is the category in which
the objects are the sets, the morphisms between objects $X$ and
$Y$ are the homotopy maps $X \times I \to Y$.

Similarly, let {\bf 2-Pre}($\mathcal{C}$) denote the 2-category in
which the objects are the presheaves, the morphisms are the maps
between presheaves, the 2-arrows are the maps $X \times I \to Y$
where $X$ and $Y$ are two presheaves and $I$ is the constant
presheaf of 2-point set $\{0, 1\}$.

Suppose that $G$ is a 2-groupoid. For any set-valued 2-functor
\cite {Mac} $X: G \to {\bf 2-Sets}$, we define its homotopy
colimit $\underrightarrow{\mbox{holim}}_G X$ as the simplicial set
with $n$-simplices:
$$\bigsqcup_{(a_0, a_1, ..., a_n)} X(a_0) \times BG(a_0, a_1)_0
\times BG(a_1, a_2)_1 \times \cdots \times BG(a_{n-1},
a_n)_{n-1}$$ where $(a_0, ..., a_n) \in \mbox{Ob}(G)^{n+1}$,
$G(a_i, a_{i+1}), i = 0, 1, ..., n-1$ are the groupoids of
morphisms from $a_i$ to $a_{i+1}$, $BG(a_i, a_{i+1})_i, i = 0, 1,
..., n-1$ is the set of all $i$-simplices in the classifying space
$BG(a_i, a_{i+1})$.

It follows that the category of 2-functor $X: G \to {\bf 2-Sets}$
and 2-natural transformations is equivalent to the category of
simplicial set maps $\pi: Y \to \overline{W}G$ such that $Y$ is
the homotopy colimit of $X$, with fibrewise maps over
$\overline{W}G$ as morphisms.

This equivalence is natural, and therefore gives an internal
description in category of presheaves of presheaf-valued 2-functor
$X: G \to {\bf 2-Pre}(\mathcal{C})$ defined on a small presheaf of
2-groupoids $G$ over a site $\mathcal{C}$.

Suppose that $G$ is a presheaf of 2-groupoids on a small site
$\mathcal{C}$. Analogous to the definition in the case of sheaf of
groupoids, when $G$ is a presheaf of 2-groupoids, we define a
\emph{G-torsor} to be a simplicial presheaf map $Y \to
\overline{W}G$ such that $Y$ is a homotopy colimit
$\underrightarrow{\mbox{holim}}_G X$ and the canonical map $Y \to
\ast$ is a local weak equivalence.

When we explicitly express the homotopy colimit
$\underrightarrow{\mbox{holim}}_G X$ by Moore complex, it's
obvious that it satisfies a pullback condition
$$\xymatrix{
\cdot\cdot\cdot\ar[r]&E \times_{G_0} G_1 \times_{G_0} G_2\ar[r]
\ar[d]&E \times_{G_0} G_1\ar[r]^-{\pi_1} \ar[d]_{\pi_2}
& E\ar[d]^f\\
\cdot\cdot\cdot\ar[r]& G_1 \times_{G_0} G_2\ar[r]& G_1\ar[r]_t &
G_0}$$ In fact, this diagram is a simplicial presheaf map $Y \to
\overline{W}G$ which satisfies the pullback condition. The map $Y
\to \ast$ is a local weak equivalence, so the presheaf $Y_0 = E$
is non-empty over $G_0$ and the action $a: E \times_{G_0} G_1 \to
E$ is transitive, $G_2$ also acts freely and transitively on $E
\times_{G_0} G_1$. Denote by $\mathbb{E}$ the 2-groupoid:
$$\xymatrix{
E \times_{G_0} G_2 \ar[r]^-s \ar@<-1ex>[r]_-t \ar[d]^-s \ar@<-1ex>[d]_-t& E \ar@/_/[l]_-i\ar@{=}[d]\\
E \times_{G_0} G_1 \ar[r]^-s \ar@<-1ex>[r]_-t\ar@/_/[u]_-i
 & E\ar@/_/[l]_-i}$$
Then $Y = \overline{W}\mathbb{E}$.

The set of $G$-torsors and natural transformations between them
forms a category. We denote it by Tors($\ast, G$). An argument
similar to that in \cite {J5} shows that the morphisms between
torsors are isomorphisms, so the category Tors($\ast, G$) is a
groupoid. Denote the path components of Tors($\ast, G$) by
$\pi_0$Tors($\ast, G$).

The set [$\ast, \overline{W}G$] of maps from $\ast$ to
$\overline{W}G$ in Ho({\bf S}Pre($ \mathcal{C}$)) may be described
as a filtered colimit by Verdier hypercovering characterization
$$[\ast, \overline{W}G] \cong \lim_{\overrightarrow{V}} \pi(V, \overline{W}G)$$
indexed over simplicial homotopy classes represented by locally
trivial fibrations (hypercovers) $V \to \ast$, where $\pi(~,~)$
indicates simplicial homotopy classes of maps.

When $G$ is a 2-groupoid, $BG$ is a simplicial groupoid and a
fibrant object in the category $s${\bf Gd} since $BG_0$ is just
the forgetful groupoid $G_0$ of the 2-groupoid $G$, the map $BG
\to \ast$ has the path lifting property; and the map $BG(x,x) \to
\ast$ is a fibration of simplicial sets since $BG(x,x) =
B(G(x,x))$ is the classifying space of the groupoid $G(x,x)$. The
functor $\overline{W}: s{\bf Gd} \to {\bf S}$ preserves fibrations
and weak equivalences (Theorem V.7.8.(2), \cite {G-J}), So
$\overline{W}G$ is a fibrant object in the category {\bf S}.

When $G$ is a presheaf of 2-groupoids, $\overline{W}G$ is a
locally fibrant object in the category {\bf S}Pre($\mathcal{C}$).

All $G$-torsors $Y \to \overline{W}G$ are of the form $
\overline{W}f: \overline{W}I \to \overline{W}G$ where $f: I \to G$
is a 2-functor defined on a 2-groupoid $I$ which is trivial in the
sense that $I(x,x) \simeq \ast$ the terminal object of the
category of presheaves of groupoids {\bf Gpd}Pre($ \mathcal{C}$)
for all local choices of objects $x$ and $I$ is locally connected.
Then $Y = \overline{W}I$ is locally fibrant, and $Y \to \ast$ is a
(local) weak equivalence. So $Y$ is a locally trivial fibrant
object, i.e., $Y$ is a hypercover of $\ast$.

Every $G$-torsor $Y \to \overline{W}G$ has within it a hypercover
$Y \to \ast$, and every morphism $Y \to Y^\prime$ of torsors is a
morphism of hypercovers.

Send the $G$-torsor $f: Y \to \overline{W}G$ to its homotopy class
$[f]: Y \to \overline{W}G$. $\ast$ is the colimit of all its
hypercovers, so $[f]$ has a unique factorization
$$\xymatrix{
Y \ar[rr]^-{[f]} \ar[dr] & &\overline{W}G \\
&\ast \ar[ur]_-{[f^\prime]}}
$$
where $[f^\prime]$ is the homotopy class of maps $\ast \to
\overline{W}G$. When we send $f$ to $[f^\prime]$ we have a map
$$ \mbox{Ob~Tors}(\ast, G) \to [\ast, \overline{W}G]$$
For any two $G$-torsors $Y \to \overline{W}G$ and $Y^\prime \to
\overline{W}G$ in same component of category Tors($\ast, G$),
there is an isomorphism $i: Y \to Y^\prime$ such that the diagrams
$$\xymatrix{
Y \ar[rr] \ar[dr]_i && \overline{W}G, & &Y \ar[rr] \ar[dr]_i && \ast \\
&Y^\prime \ar[ur]&&&&Y^\prime \ar[ur]}
$$
commute. The two torsors $f_1: Y \to \overline{W}G$ and $f_2:
Y^\prime \to \overline{W}G$ map into the same morphism
$[f^\prime]: \ast \to \overline{W}G$, so the above map can factor
through a function
$$\pi_0 \mbox{Tors}(\ast, G) \to [\ast, \overline{W}G]$$

Following Theorem 14 in \cite {J5}, we have


\begin{theorem}\label{T:2gtor}
The function $\pi_0 \mbox{Tors}(\ast, G) \to [\ast,
\overline{W}G]$ is a natural bijection for any presheaf of
2-groupoids $G$ on any Grothendieck site $\mathcal{C}$.
\end{theorem}
\begin{proof}
We have constructed the map $\varphi: \pi_0 \mbox{Tors}(\ast, G)
\to [\ast, \overline{W}G]$. Now we need find its inverse function.

For a 2-groupoid $G$ and any object $x_0 \in G$, we define the
comma 2-category $G\downarrow x_0$ to be the 2-category in which
the objects are the morphisms $x \to x_0$, the morphisms
(1-arrows) and the 2-arrows are the commutative diagrams
$$\xymatrix{
x \ar[rr] \ar[dr] && y \ar[dl]&\mbox{and}&x \ar@/^/[rr]_{id} \ar@/_/[rr] \ar[dr] & &y \ar[dl] \\
&x_0&&&&x_0}
$$
respectively, where $\xymatrix@1{x \ar@/^/[r]_{id} \ar@/_/[r] &
y}$ is the identity 2-arrow. In fact, $G\downarrow x_0$ is a
2-groupoid. Since there is only one morphism between any two
objects in $G\downarrow x_0$ and all 2-arrows are identities, the
2-groupoid $G\downarrow x_0$ is a trivial 2-groupoid, so $
\overline{W} (G\downarrow x_0)$ is contractible.

Suppose given a 2-functor $f: I \to G$, where $I$ is a trivial
2-groupoid (i.e., $I$ has only one 1-cell between any two objects
and only identity 2-cell from a 1-cell to itself). We define the
comma 2-groupoid $f \downarrow x_0$ to be the 2-groupoid in which
the objects are the pairs $(x, f(x) \to x_0), x \in I$, the
morphisms (1-arrows) and the 2-arrows are the commutative diagrams
$$\xymatrix{
x \ar[d]^{i;}  &f(x) \ar[rr]^{f(i)} \ar[dr] && f(y)
\ar[dl]&\mbox{and}& x\ar@/^/[d]_{id}^; \ar@/_/[d] &
f(x) \ar@/^/[rr]_{id} \ar@/_/[rr] \ar[dr] && f(y) \ar[dl] \\
y &&x_0&&&y&&x_0}
$$
respectively. $I$ is trivial, so is each component of $f
\downarrow x_0$. Thus $\overline{W}(f \downarrow x_0)$ is weakly
equivalent to the constant simplicial set $\pi_0\overline{W}(f
\downarrow x_0)$.

Given a hypercover $V \to \ast$, the simplicial presheaf map $V
\to \overline{W}G$ is of the form $ \overline{W}f: \overline{W}I
\to \overline{W}G$ where $f: I \to G$ is a presheaf-valued
2-functor defined on the trivial presheaf of 2-groupoids $I$.

Make the homotopy colimit construction
$$\xymatrix{
\overline{W}I(U) \ar[d]_f &d(\bigsqcup_{x_0 \to \cdots \to x_n}
\overline{W} (f\downarrow x_0) )\ar[d] \ar[l]_-\alpha
\ar[r]^-\beta & d(\bigsqcup_{x_0 \to \cdots \to x_n}
\pi_0\overline{W}
(f\downarrow x_0) )\ar[d]^{h(f)}\\
\overline{W}G(U)&d(\bigsqcup_{x_0 \to \cdots \to x_n} \overline{W}
(G(U)\downarrow x_0) ) \ar[l]_-\alpha \ar[r]^-\beta &
d(\bigsqcup_{x_0 \to \cdots \to x_n} \ast)}$$ on the presheaf
level in each section, where $x_0 \to \cdots \to x_n$ is an
$n$-simplex of simplicial set $ \overline{W}(G(U))$.

In the simplicial set $d(\bigsqcup_{x_0 \to \cdots \to x_n}
\overline{W} (f\downarrow x_0) )$, there is only one 1-simplex
between any two objects $(y, f(y) \to x)$ and $(y^\prime,
f(y^\prime) \to x^\prime)$ where $y, y^\prime \in I, x, x^\prime
\in G$
$$\xymatrix{
y \ar[d]^{i;}  &f(y) \ar[d]^{f(i)} \ar[r] & x \ar[d] \\
y^\prime&f(y^\prime)\ar[r]&x^\prime}$$ All 2-simplices are the
identities 2-arrows in $I$, and so on. So $d(\bigsqcup_{x_0 \to
\cdots \to x_n} \overline{W} (f\downarrow x_0) )$ is equivalent to
$\ast$. The top map $\alpha$ is the standard weak equivalence
associated to the simplicial map $f: \overline{W}I(U) \to
\overline{W}G(U)$, and the top map $\beta$ is also a weak
equivalence since $\overline{W} (f\downarrow x_0)$ is weak
equivalent to constant simplicial set $\pi_0\overline{W}(f
\downarrow x_0)$. So $d(\bigsqcup_{x_0 \to \cdots \to x_n}
\pi_0\overline{W} (f\downarrow x_0) )$ is equivalent to $\ast$.
$d(\bigsqcup_{x_0 \to \cdots \to x_n} \ast)$ just is the
simplicial set $ \overline{W}(G(U))$, $d(\bigsqcup_{x_0 \to \cdots
\to x_n} \pi_0\overline{W} (f\downarrow x_0) )$ is the homotopy
colimit over $G$, so the presheaf map $h(f)$ is a $G$-torsor.

When $f: \overline{W}I \to \overline{W}G$ is a $G$-torsor, $
\overline{W}I$ is the homotopy colimit of presheaves over $G$ for
some 2-functor $X$ such that $ \overline{W}I =
\underrightarrow{\mbox{holim}}_G X$ and there is a natural
isomorphism $\pi_0 \overline{W}(f\downarrow x_0) \cong X(x_0)$, so
the map $h(f)$ is canonically isomorphic to $f$.

Any homotopy $ \overline{W}I \times \triangle^1 \to \overline{W}G$
can extend to a map $ \overline{W}I \times \overline{W}G({\bf 1})
\to \overline{W}G$ where $G({\bf 1})$ is the trivial 2-groupoid on
two objects. If $I$ is locally connected and trivial, then so is
$I \times G({\bf 1})$. Corresponding to the diagram
$$\xymatrix{
\overline{W}I \ar[d]_{d^0} \ar[dr]^f\\
\overline{W}(I\times G({\bf 1}))\ar[r]^-g& \overline{W}G\\
\overline{W}I \ar[u]^{d^1} \ar[ur]_{f^\prime}}$$ there is a
diagram
$$\xymatrix{
d(\bigsqcup_{x_0 \to \cdots \to x_n} \pi_0\overline{W}(f\downarrow x_0)) \ar[d]_{d^0} \ar[dr]^{h(f)}\\
d(\bigsqcup_{x_0 \to \cdots \to x_n} \pi_0\overline{W}(g
\downarrow x_0))\ar[r]^-{h(g)}
& d(\bigsqcup_{x_0 \to \cdots \to x_n} \ast)\\
d(\bigsqcup_{x_0 \to \cdots \to x_n} \pi_0\overline{W}({f^\prime}
\downarrow x_0)) \ar[u]^{d^1} \ar[ur]_{h(f^\prime)}}$$ So any
homotopy $f \simeq f^\prime$ determines torsors $h(f)$ and
$h(f^\prime)$ which are in the same component of Tors($\ast, G$).

For any homotopy map class in the set $[\ast, \overline{W}G]$, it
has a representative $[f]: V \to \overline{W}G$. Suppose that it
has other representative $[f^\prime]: V^\prime \to \overline{W}G$.
All hypercovers over $\ast$ form a filtered category. So there
exists hypercover $V''$ and arrows $V \to V''$ and $V' \to V''$
\cite {Mac} such that the diagram
$$\xymatrix{
V \ar@{.>}[d] \ar[dr]\ar[drr]^f\\
V''\ar[r]&\ast \ar[r] &\overline{W}G \\
V' \ar@{.>}[u]\ar[ur]\ar[urr]_{f^\prime}}$$ commutes. By same
argument as above, the torsors $h(f)$ and $h(f^\prime)$ are in the
same component of Tors($\ast, G$).

There is a well-defined function
$$\psi: [\ast, \overline{W}G] \to \pi_0 \mbox{Tors}(\ast, G)$$
sending the homotopy map class $[f]$ to the $G$-torsor $h(f)$.

For any element in $ \pi_0 \mbox{Tors}(\ast, G)$, take a
representative $f: Y \to \overline{W}G$, $[f]$ is the
representative of its image in $[\ast, \overline{W}G]$, $h(f)$ is
isomorphic to $f$, so the composite function $\psi\varphi: \pi_0
\mbox{Tors}(\ast, G)\to [\ast, \overline{W}G] \to \pi_0
\mbox{Tors}(\ast, G)$ is identity.

In the above homotopy colimit diagram, the bottom $\beta$ is also
a weak equivalence because each summand $\overline{W}(G(U)
\downarrow x_0)$ is contractible. The simplicial set
$$X(U) =d(\bigsqcup_{x_0 \to \cdots \to x_n}
\overline{W} (G(U)\downarrow x_0) )$$ consists of strings $(y,x)$
of (n-)arrows
$$y_0 \to \cdots \to y_n \to x_0 \to \cdots \to x_n$$
of length $2n+1$. In fact, $y_0 \to \cdots \to y_n$ and $x_0 \to
\cdots \to x_n$ are two $n$-simplices in $\overline{W} (G(U))$.
The map $\alpha$ sends this simplex to the string $y_0 \to \cdots
\to y_n$, while $\beta$ maps it to the string $x_0 \to \cdots \to
x_n$. The $n$-simplices of the simplicial set $X(U)$ can be
identified with the simplicial maps $\triangle^n \ast \triangle^n
\to \overline{W}(G(U))$ defined on the join $\triangle^n \ast
\triangle^n$, and with simplicial structure maps induced by
precomposition with maps $\theta \ast \theta : \triangle^m \ast
\triangle^m \to \triangle^n \ast \triangle^n$. The maps $\alpha$
and $\beta$ are induced by the inclusions $\triangle^n \to
\triangle^n \ast \triangle^n$ of the left and right standard
$n$-simplex respectively.

The $n$-simplex in $X(U)$ can be expressed by diagram
$$\xymatrix{
y_0 \ar[r] \ar[d] &y_1 \ar[r] \ar[d] &\cdots \ar[r] &y_n \ar[d]\\
x_0 \ar[r] &x_1 \ar[r] & \cdots \ar[r] &x_n}$$ Then maps $\alpha$
and $\beta$ are homotopic maps. This construction is natural,
hence the two maps represent the same morphism in the homotopy
category of simplicial presheaves.

So the two composition maps
$$\xymatrix@1{d(\bigsqcup_{x_0 \to
\cdots \to x_n} \overline{W} (f \downarrow x_0)
)\ar[r]&d(\bigsqcup_{x_0 \to \cdots \to x_n} \overline{W}
(G(U)\downarrow x_0) )\ar[r]^-\alpha & d(\bigsqcup_{x_0 \to \cdots
\to x_n}\ast)}$$ and $\xymatrix@1{d(\bigsqcup_{x_0 \to \cdots \to
x_n} \overline{W} (f \downarrow x_0) )\ar[r]&d(\bigsqcup_{x_0 \to
\cdots \to x_n} \overline{W} (G(U)\downarrow x_0) )\ar[r]^-\beta &
d(\bigsqcup_{x_0 \to \cdots \to x_n}\ast)}$ are homotopic. Denote
this homotopy class as $[g]$.

$d(\bigsqcup_{x_0 \to \cdots \to x_n} \overline{W} (f \downarrow
x_0) )$ can be written as $\overline{W}H$ for some 2-groupoid $H$,
hence it's fibrant, and it's weak equivalent to $\ast$, then the
simplicial presheaf $d(\bigsqcup_{x_0 \to \cdots \to x_n}
\overline{W} (f \downarrow x_0) )$ is a hypercover over $\ast$. By
the above homotopy colimit diagram, $[f], [g]$ and $[h(f)]$
represent the same element in $[\ast, \overline{W}G]$.

The composition function $\varphi\psi: [\ast, \overline{W}G] \to
\pi_0 \mbox{Tors}(\ast, G) \to [\ast, \overline{W}G] $ sends $[f]$
to $[h(f)]$, so it's identity.
\end{proof}


\subsection{Torsors for presheaves of simplicial groups}

Suppose that $G$ is a presheaf of simplicial groups on a small
site $\mathcal{C}$.

Let $X$ be a simplicial presheaf. $G$ \emph{acts} on $X$ if there
is a morphism of simplicial presheaves
$$\mu: G \times X \to X$$
so that the following diagrams commute:
$$\xymatrix{
G \times G \times X \ar[r]^-{1 \times \mu} \ar[d]_-{m \times 1} &
G \times X \ar[d]^-\mu\\
G \times X \ar[r]_-\mu &X}
$$
and
$$\xymatrix{
X \ar[d]_-i \ar[dr]^-{1_X}\\
G \times X \ar[r]_-\mu & X}
$$
where $m$ is the multiplication in $G$ and $i(X) = (e, X)$ on each
level of each section.

Such a simplicial presheaf $X$ is called a \emph{simplicial
$G$-presheaf}. Let {\bf S}Pre($ \mathcal{C})_G$ be the category of
simplicial $G$-presheaves.

The forgetful functor {\bf S}Pre($ \mathcal{C})_G \to $ {\bf
S}Pre($ \mathcal{C})$ has a left adjoint given by
$$ X \mapsto G \times X$$
It's easy to prove that there is a closed model structure on the
category {\bf S}Pre($ \mathcal{C})_G$, where a map $f: X \to Y$ of
simplicial $G$-presheaves is a fibration (respectively weak
equivalence) if the underlying map of simplicial presheaves is a
global fibration (respectively local weak equivalence) and a
cofibration is a map which has the left lifting property with
respect to all both fibrations and weak equivalences (cf. \cite
[Theorem V.2.3.]{G-J}).

A \emph{principal $G$-bundle} (or \emph{principal $G$-fibration})
$f: E \to B$ is a local fibration in {\bf S}Pre($ \mathcal{C})_G$
so that
\begin{itemize}
\item[(1)] $B$ has trivial $G$-action;
\item[(2)] $E$ is a cofibrant simplicial $G$-presheaf, and
\item[(3)] the induced map $E/G \to B$ is an isomorphism.
\end{itemize}

\begin{lemma}
Suppose that $X$ is a cofibrant simplicial $G$-presheaf. Then $G$
acts freely on $X$ in all sections.
\end{lemma}
\begin{proof}
Suppose that the functor $L_U: {\bf S} \to {\bf S}Pre(
\mathcal{C})$ is the functor $?_U$ in \cite {J1}, the left adjoint
functor of the $U$-section functor $X \to X(U)$:
$$L_U(Y)(V) = \bigsqcup_{\varphi: V \to U} Y.$$
The maps $G \times L_U\partial\triangle^n \to G \times
L_U\triangle^n$ generate the cofibrations of the category of
simplicial $G$-presheaves, and any pushout
$$\xymatrix{
G \times L_U\partial\triangle^n \ar[r]\ar[d] &Z\ar[d]\\
G \times L_U\triangle^n \ar[r] &W}
$$
has the effect of adding some freely generated $G(U)$-space to
$Z(U)$ for each $U \in \mathcal{C}$. The cofibration $\emptyset
\to X$ has a factorization
$$\xymatrix{
\emptyset \ar[r] \ar[dr] & V\ar[d]^-\pi\\
 & X}
$$
where $\pi$ is a trivial fibration and the map $\emptyset \to V$
is a transfinite colimit of pushouts of the above form. It follows
that $G$ acts freely on $V$ in all sections. But then, by a
standard argument $X$ is a retract of $V$ (since $\pi$ is a
trivial fibration), so that $G$ acts freely on $X$ in all
sections.
\end{proof}

Now for the converse:

\begin{lemma}
Suppose that $G$ acts freely on the simplicial $G$-presheaf $Y$ in
all sections. Then $Y$ is cofibrant.
\end{lemma}
\begin{proof}
Consider the partially ordered set of all cofibrant subobjects $K
\subset Y$, where $K \leq L$ if there is a $G$-cofibration $K
\subset L$ which respects the inclusions into $Y$. This poset is
non-empty since $G\langle x\rangle \subset Y$ is a cofibrant
subobject of $Y$ by the condition that $G$ acts freely on $Y$ in
all sections for any simplex $x \in Y(U)$. If
$$K_1 \leq K_2 \leq \cdots$$
is a totally ordered collection of cofibrant subobjects, the
$K_\infty = \cup K_i$ is cofibrant and all inclusions $K_i \subset
K_\infty$ are cofibrations. It follows that the poset of cofibrant
subobjects of $Y$ is inductively ordered.

Zorn's lemma therefore asserts that the poset of cofibrant
subobjects of $Y$ has maximal elements. Pick such a maximal
subobject $M$ and assume that $M \neq Y$. Then there is a simplex
$x \in Y(U) - M(U)$ of minimal dimension, and the diagram
$$\xymatrix{
G(\langle x \rangle) \cap M \ar[r]\ar[d] & M \ar[d]\\
G(\langle x \rangle) \ar[r] & G(\langle x \rangle) \cup M}
$$
is a pushout of $G$-subcomplexes of $Y$. But
$$G(\langle \partial x \rangle) = G(\langle x \rangle) \cap M$$
where $\langle \partial x\rangle$ is the subcomplex generated by
the faces of $x$. It follows that the inclusion map $M \to
G(\langle x \rangle) \cup M$ is a cofibration of simplicial
$G$-presheaves, which contradicts the maximality of $M$ if $M \neq
Y$. Thus, $Y$ is cofibrant.
\end{proof}

Hence a simplicial $G$-presheaf $X$ is cofibrant if and only if
$G$ acts freely on it in all sections.

For every cofibrant simplicial $G$-presheaf $X$ the canonical map
$X(U) \to X(U)/G(U)$ is a fibration of simplicial sets \cite
[Corollary V.2.7] {G-J} for any $U \in \mathcal{C}$. Then $X \to
X/G$ is a local fibration \cite [Corollary 1.8] {J1} where $X/G$
is the simplicial presheaf $U \mapsto X(U)/G(U)$. So the canonical
map $X \to X/G$ is a principal $G$-bundle. Any $G$-bundle $f: E
\to B$ is isomorphic to a quotient map
$$ q: X \to X/G$$
where $X \in {\bf S}Pre(\mathcal{C})_G$ is cofibrant. That means
that cofibrant simplicial $G$-presheaves $X$ correspond to
principal $G$-bundles $X \to X/G$.

A \emph{$G$-torsor} on $\mathcal{C}$ is a cofibrant simplicial
$G$-presheaf $X$ such that the canonical map $X/G \to \ast$ is a
hypercover. A map $f: X \to Y$ of $G$-torsors is just a
$G$-equivariant map of simplicial presheaves. Write {\bf G-tors}
for the corresponding category and $\pi_0({\bf G-tors})$ for the
corresponding path components of {\bf G-tors}.

Choose a factorization
$$\xymatrix{
\emptyset \ar[r]^-i \ar[dr] & EG \ar[d]^-\pi\\
& \ast}
$$
where $i$ is a cofibration and $\pi$ is a trivial fibration in the
category of simplicial $G$-presheaves {\bf S}Pre($\mathcal{C})_G$.
Write $BG = EG/G$. Observe that $BG$ is locally fibrant, since it
is the presheaf of Kan complexes (cf. Lemma V.3.7. \cite {G-J}).

Note that $EG$ is unique up to equivariant homotopy equivalence,
so $q: EG \to BG$ is unique up to homotopy equivalence.

Suppose that $U \to BG$ is a map of simplicial presheaves, where
$U \to \ast$ is a hypercover. Then the pullback $U \times_{BG} EG$
has a free $G$-action and is therefore a cofibrant simplicial
$G$-presheaf. The projection map $U \times_{BG} EG \to U$ is the
quotient by the $G$-action and $U \to \ast$ is a hypercover, so
that $U \times_{BG} EG$ is a $G$-torsor.

Further, any string of morphisms
$$U_1 \to U_2 \to BG$$
where $U_1 \to U_2$ is a map of hypercovers induces a morphism
$$U_1\times_{BG} EG \to U_2 \times_{BG} EG$$
of $G$-torsors in the obvious way.

\begin{remark}
In the standard theory $WG$ is a cofibrant simplicial
$G$-presheaf, but the map $WG \to \ast$ is only a (trivial) local
fibration. Find a factorization
$$\xymatrix{
WG \ar[r]^-j \ar[dr] &EG\ar[d]^-\pi\\
& \ast}
$$
such that $\pi$ is a trivial fibration and $j$ is a trivial
cofibration in the category of simplicial $G$-presheaves {\bf
S}Pre($\mathcal{C})_G$. Then the induced comparison of principal
$G$-bundles
$$\xymatrix{
WG \ar[d]\ar[r]^-j_-\simeq &EG \ar[d]\\
\overline{W}G \ar[r]_-{j_\ast} &BG}
$$
induces a local weak equivalence $j_\ast: \overline{W}G \to BG$
such that $j_\ast$ has a homotopy inverse $j^\ast: BG \to
\overline{W}G$ where $j^\ast$ is a local weak equivalence as well.
Similarly, if $U \to \overline{W}G$ is a map of simplicial
presheaves, where $U \to \ast$ is a hypercover, then $U
\times_{\overline{W}G} WG$ is also a $G$-torsor.

A similar argument works for the diagonal map $d(EG) \to d(BG)$
induced by the standard bisimplicial presheaf map $EG \to BG$,
because the diagonal map is a principal $G$-fibration and the
object $d(EG)$ is weak equivalence to a point. It follows that
$d(BG) \simeq BG$ for the two different senses of $BG$.
\end{remark}

A simplicial presheaf $X$ is called \emph{projective fibrant} if
the map $X \to \ast$ has the right lifting property with respect
to all maps $L_U \wedge^n_k \to L_U \triangle^n, U \in
\mathcal{C}$. It's obvious that the simplicial presheaf $BG$ is
projective fibrant. In effect, the map $EG \to BG$ is a principal
$G$-bundle and hence a surjective Kan fibration in each section,
and $EG$ is globally fibrant, and hence a Kan complex in each
section, so that $BG$ is a Kan complex in each section. The
simplicial presheaf $\overline{W}G$ is projective fibrant as well.

Say that a presheaf of simplicial groupoids $H$ is projective
fibrant if the object $\overline{W}H$ is projective fibrant.

More generally, one says that the map $p: G \to H$ is a projective
fibration of presheaves of simplicial groupoids if the induced map
$\overline{W}G \to \overline{W}H$ is a projective fibration. This
is equivalent to the assertion that the map $p$ has the right
lifting property with respect to all maps $G(L_U \wedge^n_k) \to
G(L_U \triangle^n), U \in \mathcal{C}$.

Every presheaf of simplicial groupoids $H$ has a \emph{projective
fibrant model} $i: H \to H_f$, in the sense that the map $i$ is a
sectionwise weak equivalence which has the left lifting property
with respect to all projective fibrations, and $H_f$ is projective
fibrant. This is a consequence of the obvious small object
argument.

Suppose once again that $G$ is a presheaf of simplicial groups.
Then $G$ is projective fibrant by a standard argument. Suppose
that given a homotopy $h: U \times \triangle^1 \to \overline{W}G$
where $U$ is a hypercover. Then the induced map $h_\ast: G(U
\times \triangle^1) \to G$ has a factorization
$$
\xymatrix{G(U \times \triangle^1) \ar[d]_-j\ar[r]^-{h_\ast}& G\\
    H \ar[ur]_-\pi}
$$
where the map $j: G(U \times \triangle^1) \to H$ is a projective
fibrant model for $G(U \times \triangle^1)$ in the category of
presheaves of simplicial groupoids. $j$ has the left lifting
property with respect to all projective fibrations, $G$ is
projective fibrant, so the lifting $\pi$ exists. It follows that $
\overline{W}H \to \ast$ is a hypercover since it's sectionwise
weak equivalence and fibration, hence it's a locally trivial
fibration, and there is a commutative diagram
$$
\xymatrix{U\ar[d]\ar[dr]^-{hd^0}\\
\overline{W}H \ar[r] & \overline{W}G\\
U\ar[u]\ar[ur]_-{hd^1}}
$$
In particular, the principal $G$-fibrations over $U$ which are
induced by the maps $hd^1$ and $hd^0$ are in the same component of
the $G$-torsor category.

For the composite map $\xymatrix@1{U \ar[r]^-f & BG
\ar[r]^-{j^\ast} & \overline{W}G}$, the principal $G$-bundles $U
\times_{BG} EG$ and $U \times_{\overline{W}G} WG$, which are
induced by the maps $f$ and $j^\ast f$ respectively, are in the
same component of the $G$-torsor category. If $f, g: U \to BG$ are
homotopic, so are $j^\ast f, j^\ast g: U \to \overline{W}G$. Hence
the $G$-torsors which are induced by $f$ and $g$ are in same
component of {\bf G-tors}.

It follows that there is a well defined function
$$\psi_G: [\ast, BG] = \lim_{\longrightarrow \atop U}\pi(U,
BG) \to \pi_0({\bf G-tors})$$ which is given by sending the naive
homotopy class of a map $U \to BG$ defined on a hypercover $U \to
\ast$ to the path component of the object $U \times_{BG} EG$.

Suppose that $X$ is a $G$-torsor. Then $X$ is a cofibrant
simplicial $G$-presheaf so that the lifting $\phi$ exists in the
diagram
$$\xymatrix{
\emptyset \ar[r]\ar[d] &EG \ar[d]^-\pi\\
X \ar[r]\ar[ur]^-\phi & \ast}
$$
since $\pi$ is a trivial fibration. Moreover, any two such
liftings are homotopic. Make a fixed choice of lifting $\phi_X$
for all $G$-torsors $X$, and write $\phi_{X_\ast}: X/G \to BG$ for
the corresponding induced map.

Suppose that $f: X \to Y$ is a morphism of $G$-torsors. Then
$\phi_Y f$ is a second possible choice for $\phi_X$ and so $\phi_Y
f$ and $\phi_X$ are (naively) equivariantly homotopic. It follows
that the maps $\phi_{X_\ast}$ and $\phi_{Y_\ast} f_\ast$ are
naively homotopic. The assignment
$$X \mapsto [\phi_{X_\ast}] \in \pi(X/G, BG)$$
therefore determines a well-defined function
$$\varphi_G: \pi_0({\bf G-tors}) \to \lim_{\longrightarrow \atop {U \to \ast}}\pi(U,
BG).$$ It's obvious that both maps $\varphi_G \cdot \psi_G$ and
$\psi_G \cdot \varphi_G$ are identity maps.

We have therefore shown

\begin{theorem}\label{T:Gsheaf}
The function $\psi_G: [\ast, BG] \to \pi_0({\bf G-tors})$ is a
bijection for all presheaves of simplicial groups $G$ on a
Grothendieck site $\mathcal{C}$.
\end{theorem}

Note that the theorem is independent of the choice of the
cofibrant model $EG$ for the point.

\begin{corollary}\label{C:Gsheaf}
Suppose that $G$ is a presheaf of simplicial groups on a
Grothendieck site $\mathcal{C}$. Then there is a bijection
$$[\ast, \overline{W}G] \cong \pi_0({\bf
G-tors})$$
\end{corollary}

\begin{remark}
Every trivial cofibration $i: A \to B$ of simplicial
$G$-presheaves induces a trivial cofibration $i_\ast: A/G \to
B/G$. In effect, $i$ has the left lifting property with respect to
all global fibration $p: X \to Y$ of simplicial presheaves with
trivial $G$-action.

Suppose that $X$ is a cofibrant simplicial $G$-presheaf such that
the induced map $X/G \to \ast$ is a weak equivalence. Find a
trivial cofibration $j: X \to \tilde{X}$ in the category of
simplicial $G$-presheaves such that $\tilde{X}$ is fibrant. Then
the induced map $j_\ast: X/G \to \tilde{X}/G$ is a trivial
cofibration of simplicial $G$-presheaves, and $\tilde{X}/G$ is
locally fibrant so the map $\tilde{X}/G \to \ast$ is a hypercover.
Write {$\bf G-tors_0$} for the category of cofibrant simplicial
$G$-presheaves $X$ such that $X/G \to \ast$ is a local weak
equivalence. Then the inclusion
$$ {\bf G-tors}\subset {\bf G-tors_0}$$
induces an isomorphism
$$ \pi_0({\bf G-tors})\cong \pi_0({\bf G-tors_0}).$$
\end{remark}

\begin{remark}
Write {$\bf G-tors_1$} for the category of simplicial
$G$-presheaves $Y$ such that the canonical map $d(EG \times_G Y)
\to \ast$ is a local weak equivalence, where $d(X)$ denotes the
diagonal of a bisimplicial object $X$. Then there is an inclusion
$$ {\bf G-tors_0}\subset {\bf G-tors_1}$$
since the canonical map $d(EG \times_G Z) \to Z/G$ is a weak
equivalence if $Z$ is cofibrant. On the other hand, if $X$ is a
simplicial $G$-presheaf such that $d(EG \times_G X) \to \ast$ is a
weak equivalence, there is a trivial fibration $Z \to X$ of
simplicial $G$-presheaves such that $Z$ is cofibrant. The induced
map $d(EG \times_G Z) \to d(EG \times_G X)$ is a local weak
equivalence, so that $Z$ is an object of {$\bf G-tors_0$}. It
follows that there is an isomorphism
$$ \pi_0({\bf G-tors_0})\cong \pi_0({\bf G-tors_1})$$
If one defines \emph{G-torsors} on $\mathcal{C}$ as simplicial
presheaves $X$ with principal $G$-action such that the canonical
map $d(EG \times_G X) \simeq X/G \to \ast$ is a weak equivalence,
the similar bijection of Theorem \ref{T:Gsheaf} exists.
\end{remark}


\subsection{Torsors for presheaves of simplicial groupoids}

First of all, we extend the results in Section 3.4.

Write ${\bf Triv}/X$ for the category whose objects are all
simplicial presheaf morphisms $W \to X$ such that the map $W \to
\ast$ is a local weak equivalence. Observe that there is a
function
$$\psi_X : \pi_0({\bf Triv}/X) \to [\ast, X]$$
which is defined by associating to an object $W \to X$ the
composite
$$\xymatrix@1{\ast & W \ar[l]_-\simeq\ar[r] &X}$$
in the homotopy category.

Note that there are corresponding constructions for any object $X$
of a model category {\bf M}.

\begin{lemma}\label{L:Triv}
Suppose that {\bf M} is a right proper model category, and suppose
that the map $f: X \to Y$ is a weak equivalence. Then the induced
map
$$ f_\ast: \pi_0({\bf Triv}/X) \to \pi_0({\bf Triv}/Y)$$
is a bijection.
\end{lemma}
\begin{proof}
The function $f_\ast$ is induced by a functor which is defined by
associating to the object $W \to X$ the composite
$$\xymatrix@1{ W \ar[r] & X \ar[r]^-f &Y.}$$
Suppose that $U \to Y$ is an object of ${\bf Triv}/Y$. Choose a
factorization
$$\xymatrix{
U \ar[r]^-j \ar[dr] &V \ar[d]^-p\\
& Y}
$$
where $j$ is a trivial cofibration and $p$ is a fibration. Form
the pullback
$$\xymatrix{
X \times_Y V \ar[r]^-{f'} \ar[d] &V \ar[d]^-p\\
X \ar[r]_-f &Y}
$$
Then the map $f'$ is a weak equivalence by the properness
assumption, so that the projection $X \times_Y V \to X$ is an
object of ${\bf Triv}/X$. Observe that the path component of this
object is independent of the choice made, and is independent of
the choice of representative for the path component of $U \to Y$.
It follows that there is a well-defined function
$$ g: \pi_0({\bf Triv}/Y) \to \pi_0({\bf Triv}/X).$$
The composite functions $g\cdot f_\ast$ and $f_\ast \cdot g$ are
both identities.
\end{proof}

\begin{lemma}\label{L:Triv2}
Suppose that $Y$ is an object of a right proper model category
{\bf M} in which the terminal object $\ast$ is cofibrant. Then the
function
$$ \psi_Y: \pi_0({\bf Triv}/Y) \to [\ast, Y]$$
is a bijection.
\end{lemma}
\begin{proof}
Suppose first of all that $Y$ is fibrant. Then the function
$$\pi(\ast, Y) \to [\ast, Y]$$
is a bijection since $\ast$ is cofibrant. Here, $\pi(\ast, Y)$
denotes homotopy classes of maps with respect to a fixed cylinder
object $I$ of $\ast$. If two maps $f,g: \ast \to Y$ are homotopic,
then there is a diagram
$$\xymatrix{
\ast\ar[d]_-{d_0} \ar[dr]^-f\\
I \ar[r] &Y\\
\ast\ar[u]^-{d_1}\ar[ur]_-g}
$$
Then the morphisms $d_0$ and $d_1$ are weak equivalences, so that
$f$ and $g$ are in the same path component of ${\bf Triv}/Y$. It
follows that there is a well defined function
$$\phi: \pi(\ast, Y) \to \pi_0({\bf Triv}/Y)$$
and the diagram
$$\xymatrix{
\pi(\ast, Y) \ar[r]^-\cong \ar[dr]_-\phi &[\ast, Y]\\
&\pi_0({\bf Triv}/Y) \ar[u]_-{\psi_Y}}
$$
commutes. Finally, if $U \to Y$ is an object of ${\bf Triv}/Y$,
there is a factorization
$$\xymatrix{
U \ar[r]^-j \ar[dr] &V\ar[d]^p\\
& \ast}
$$
where $j$ is a trivial cofibration and $p$ is a trivial fibration.
The fibration $p$ has a section $s: \ast \to V$ since $\ast$ is
cofibrant, and the map $U \to Y$ extends to a map $V \to Y$ since
$j$ is a trivial cofibration and $Y$ is fibrant. It follows that
the function $\phi$ is surjective, and is therefore a bijection.

The map $\psi_Y$ is therefore a bijection if $Y$ is fibrant. The
general case follows from Lemma \ref{L:Triv}.
\end{proof}

\begin{lemma}
Suppose that $G$ is a presheaf of simplicial groups. Then there is
a bijection
$$[\ast, BG] \cong \pi_0({\bf G-tors_0}).$$
\end{lemma}
\begin{proof}
We establish the existence of a bijection
$$\pi_0({\bf Triv}/BG) \cong \pi_0({\bf G-tors_0}).$$

First of all, there is a functor ${\bf Triv}/BG \to {\bf
G-tors_0}$ which is defined by associating the cofibrant
simplicial $G$-presheaf $X \times_{BG} EG$ to the object $X \to
BG$ of ${\bf Triv}/BG$.

Suppose that $Z$ is a cofibrant simplicial $G$-presheaf such that
$Z/G \to \ast$ is a local weak equivalence. Then there is a
$G$-equivariant map $Z \to EG$ and an induced map $Z/G \to BG$.
The class of the object $Z/G \to BG$ in $\pi_0({\bf Triv}/BG)$ is
independent of the choices that have been made: any two
$G$-equivariant maps $Z \to EG$ are naively homotopic and so the
induced maps $Z/G \to BG$ are naively homotopic and hence
represent the same path component of $\pi_0({\bf Triv}/BG)$. It
follows that there is a well-defined function
$$\pi_0({\bf G-tors_0}) \to \pi_0({\bf Triv}/BG)$$
and this function is the inverse of the function in $\pi_0$ which
is induced by the functor of the previous paragraph.
\end{proof}

\begin{corollary}
There is a bijection
$$[\ast, BG] \cong \pi_0({\bf G-tors_1}).$$
\end{corollary}

Recall that the objects of the category {$\bf G-tors_1$} are
simplicial $G$-presheaves $Z$ such that $d(EG \times_G Z) \to
\ast$ is a local weak equivalence.

\begin{remark}
If $G$ is a sheaf of groups, then a $G$-torsor $X$ is naturally a
member of ${\bf G-tors_0}$ after identification of $X$ with a
constant simplicial $G$-sheaf, and in this way the category {\bf
G-tors} of ordinary $G$-torsors imbeds in ${\bf G-tors_0}$, and
the induced function
$$\pi_0({\bf G-tors}) \to \pi_0({\bf G-tors_0})~~~~~~~~~~~~~~~~~~~~~~~~~~~~~(1)$$
is a bijection.

Note that $EG$ is a cofibrant simplicial $G$-sheaf such that $EG
\to \ast$ is a global fibration and $BG = EG/G$. Write $NG$ for
the nerve of $G$ in this instance. Recall that there is a weak
equivalence $NG \to BG$.

Write $C(U)$ for the $\check{\mbox{C}}$ech resolution associated
to a covering $U \to \ast$. The objects of the category ${\bf
cov}/NG$ are the morphisms $C(U) \to NG$, and the morphisms are
just commutative diagrams. There is an obvious inclusion functor
${\bf cov}/NG \subset {\bf Triv}/NG$. Observe that there are
bijections
$$\pi_0({\bf cov}/NG) \cong \pi_0({\bf Triv}/NG) \cong \pi_0({\bf
Triv}/BG)~~~~~~~~(2)$$

Now it's well known that the set of naive homotopy classes
$\pi(C(U), NG)$ is isomorphic to the set of isomorphism classes of
$G$-torsors which trivialize over $U$, and that the set of
isomorphism classes of $G$-torsors is isomorphic to
$$\lim_{\longrightarrow \atop U}\pi(C(U), NG),$$
where the colimit is indexed over the coverings $U \to \ast$. If
two maps $C(U) \to NG$ are homotopic, the homotopy $C(U) \times
\triangle^1 \to NG$ factors through the nerve $N\pi(C(U) \times
\triangle^1)$, which is itself a $\check{\mbox{C}}$ech resolution
$C(V)$ for some covering $V \to \ast$. It follows that if two maps
$C(U) \to NG$ are homotopic, then they represent the same element
of $\pi_0({\bf cov}/NG)$, and there is a well-defined function
$$\lim_{\longrightarrow \atop U}\pi(C(U), NG) \to \pi_0({\bf cov}/NG).$$
This function is a bijection, with an obvious inverse. The
function (1) is therefore isomorphic to the composite isomorphism
(2).
\end{remark}

The category $s${\bf Gd}Pre($\mathcal{C}$) (or $s${\bf
Gd}Shv($\mathcal{C}$)) of (pre)sheaves of simplicial groupoids is
right proper (Theorem~\ref{T:simpresgpro}). The terminal object
$\ast$ is cofibrant in $s${\bf Gd}Pre($\mathcal{C}$) (or $s${\bf
Gd}Shv($\mathcal{C}$)).

\begin{lemma}\label{L:Diagram}
There is a natural diagram of bijections
$$\xymatrix{
\pi_0({\bf Triv}/G) \ar[d]_-{\overline{W}}\ar[r]^-{dB}
&\pi_0({\bf Triv}/dBG)\ar[dl]^-{j_\ast}\\
\pi_0({\bf Triv}/\overline{W}G)}
$$
for all objects $G \in s{\bf Gd}Pre(\mathcal{C})$.
\end{lemma}
\begin{proof}
The proof is essentially trivial, and follows from the existence
of the diagram of bijections
$$\xymatrix{
[\ast, G] \ar[d]_-{\overline{W}}\ar[r]^-{dB}
&[\ast, dBG]\ar[dl]^-{j_\ast}\\
[\ast, \overline{W}G]}
$$
together with Lemma \ref{L:Triv2}.
\end{proof}

A (\emph{weakly}) \emph{simplicial category} $ \mathcal{A}$ is a
simplicial object in the category of categories having a discrete
simplicial class of objects; in other words, a (weakly) simplicial
category $\mathcal{A}$ is a category enriched in simplicial sets.
The simplicial groupoid is a weakly simplicial category. The full
simplicial set category {\bf S} with the function complexes {\bf
Hom}$(X, Y)$ is also a weakly simplicial category.

The simplicial set of morphisms from $A$ to $B$ in a weakly
simplicial category $ \mathcal{A}$ is denoted by
$\mathcal{A}(A,B)$; the corresponding set of $n$-simplices
$\mathcal{A}(A,B)_n$ is the set of morphisms from $A$ to $B$ in
the category at the level $n$.

A \emph{simplicial functor} $f: \mathcal{A} \to \mathcal{B}$ is a
morphism of (weakly) simplicial categories. This means that $f$
consists of a function $f: \mbox{Ob}( \mathcal{A}) \to \mbox{Ob}(
\mathcal{B})$ and simplicial set maps $f: \mathcal{A}(A,B) \to
\mathcal{B}(f(A),f(B))$ which respect identities and compositions
at all levels.

A \emph{natural transformation} $\eta: f \to g$ of simplicial
functors $f,g: \mathcal{A} \to \mathcal{B}$ consists of morphisms
$$\eta_A: f(A) \to g(A)$$
in hom$(f(A),g(A)) = \mathcal{B}(f(A),g(A))_0$, one for each
object $A$ of $ \mathcal{A}$, such that the following diagram of
simplicial set maps commutes
$$
\xymatrix{\mathcal{A}(A,B)\ar[r]^-f\ar[d]_-g&\mathcal{B}(f(A),f(B))\ar[d]^-{\eta_{B_\ast}}\\
\mathcal{B}(g(A),g(B))\ar[r]_-{\eta^\ast_A}&\mathcal{B}(f(A),g(B))}
$$
for each pair of objects $A,B$ of $ \mathcal{A}$.

Suppose that $C$ is a simplicial category, the simplicial functor
taking values in simplicial sets $X: C \to {\bf S}$ gives rise to
a bisimplicial set $BE_C X$ with simplicial set
$$\bigsqcup_{(a_0,a_1,\cdots,a_n)} X(a_0) \times C(a_0,a_1) \times
\cdots \times C(a_{n-1},a_n)$$ in horizontal degree $n$. In
vertical degree $m$, it is the simplicial set
$\underrightarrow{\mbox{holim}}_{C_m} X_m$.

The \emph{homotopy colimit} of the simplicial functor $X$ is the
diagonal $d(BE_C X)$; one usually writes
$\underrightarrow{\mbox{holim}}_C X = d(BE_C X)$.

$E_C X$ is a translation simplicial category. In effect, each of
$m$-simplex simplicial functors $X_m : C_m \to {\bf Sets}$ gives
rise to a translation category $E_{C_m} X_m$ having objects
$(i,x)$ with $i$ an object of $C_m$ (or $C$) and $x\in X_m(i)$,
and with morphisms $\alpha: (i,x) \to (j,y)$ where $\alpha: i \to
j $ is a morphism of $C_m$ such that $X_m(\alpha)(x) = y$. Then
the nerve $BE_{C_m} X_m$ is the homotopy colimit
$\underrightarrow{\mbox{holim}}_{C_m} X_m$. Furthermore, the data
is simplicial in $n$, so the simplicial object $BE_C X$ is indeed
a bisimplicial set.

Suppose that $G$ is a simplicial group and that $X$ is a
simplicial set admitting a left $G$-action. Then the functor $X: G
\to {\bf S}$ sending the unique object $\ast$ to the simplicial
set $X$ is a simplicial functor. The Borel construction $d(EG
\times_G X)$ is the homotopy colimit
$\underrightarrow{\mbox{holim}}_G X$.

When $G$ is a presheaf of simplicial groups on a small site
$\mathcal{C}$, define $G$-torsors as the simplicial functors
taking values in simplicial presheaves $X: G \to {\bf
S}Pre{(\mathcal{C})}$ such that $\underrightarrow{\mbox{holim}}_G
X \to \ast$ is a locally weak equivalence, then the category {\bf
G-tors} of such $G$-torsors coincides with the category {$\bf
G-tors_1$} in Section 3.4.

\begin{lemma}
Suppose that $C$ is a category enriched in simplicial sets and
that $X: C \to {\bf S}$ is a simplicial functor taking values in
simplicial sets. Suppose that all arrows $a \to b$ of $C_0$ induce
weak equivalences $X(a) \to X(b)$. Then the map $X(a) \to F_a$
taking values in the homotopy fibre over $a$ of the diagonal
simplicial set map $\underrightarrow{\mbox{holim}}_C X \to d(BC)$
is a weak equivalence.
\end{lemma}
\begin{proof}
For the pullback diagram of bisimplicial sets:
$$
\xymatrix{ \underset
{a\stackrel{1}{\longrightarrow}a\stackrel{1}{\longrightarrow}
\cdots\stackrel{1}{\longrightarrow}a\in BC_n}{\bigsqcup}
 X(a)\ar[d]\ar[r]
&\underset{(a_0,a_1,\cdots,a_n)}{\bigsqcup} X(a_0)
\times C(a_0,a_1) \times \cdots \times C(a_{n-1},a_n)\ar[d]^-{\pi}\\
\ast\ar[r]_-a&BC = \underset {(a_0,a_1,\cdots,a_n)}{\bigsqcup}
C(a_0,a_1) \times \cdots \times C(a_{n-1},a_n)}
$$
applying the diagonal functor $d$ one obtains the pullback diagram
(since $d$ has a left adjoint $d^\ast$ \cite [p. 220] {G-J} then
it preserves pullback):
$$
\xymatrix{X(a) \ar[r]\ar[d]&\underrightarrow{\mbox{holim}}_C
X\ar[d]^-{d(\pi)}\\
\ast\ar[r]_-a& d(BC)}
$$

Consider all bisimplices $\sigma: \triangle^{r,s} \to BC$ of $BC$,
form the pullback diagram
$$
\xymatrix{ {\pi}^{-1}(\sigma)\ar[d]\ar[r]
&\underset{(a_0,a_1,\cdots,a_n)}{\bigsqcup} X(a_0)
\times C(a_0,a_1) \times \cdots \times C(a_{n-1},a_n)\ar[d]^-{\pi}\\
\triangle^{r,s}\ar[r]_-\sigma&BC }
$$
in the category of bisimplicial sets. Since $\triangle^{0,0} =
\ast$, then $X(a) = d\pi^{-1}(a)$.

The bisimplices $\triangle^{m,n}\to BC$ of $BC$ are the objects of
the \emph{category of bisimplices} of $BC$, denoted by
$\Delta^{\times 2} \downarrow BC$. A morphism $\sigma \to \tau$ of
this category is a commutative diagram of bisimplicial set maps
$$
\xymatrix{ \triangle^{r,s}\ar[dr]^-\sigma\ar[dd]\\
&BC\\
\triangle^{m,n}\ar[ur]_-\tau}
$$
The assignment $\sigma \mapsto \pi^{-1}(\sigma)$ defines a functor
$$\pi^{-1}: \Delta^{\times 2} \downarrow BC \to {\bf S}^2.$$

This Lemma follows from Lemma IV.5.7 in \cite {G-J} if we can show
that each morphism of bisimplices
$$
\xymatrix{ \triangle^{r,s}\ar[dr]^-{\tau}\ar[dd]_{(\zeta_1,\zeta_2)}\\
&BC\\
\triangle^{k,l}\ar[ur]_-{\sigma}}
$$
induces a weak equivalence
$$(\zeta_1,\zeta_2)_\ast: d\pi^{-1}(\tau) \to d\pi^{-1}(\sigma)$$
Since the argument in p. 246 of \cite {G-J} shows that if the
pullback diagram
$$
\xymatrix{ d\pi^{-1}(a)
\ar[r]\ar[d]&\underrightarrow{\mbox{holim}}_{\triangle^{\times 2}
\downarrow BC} d\pi^{-1}\ar[d]^-{\pi}\\
\ast\ar[r]_-a& B(\triangle^{\times 2} \downarrow BC)}
$$
of simplicial sets is homotopy cartesian, then the pullback
diagram
$$
\xymatrix{X(a) \ar[r]\ar[d]&\underrightarrow{\mbox{holim}}_C
X\ar[d]^-{d(\pi)}\\
\ast\ar[r]_-a& d(BC)}
$$
is homotopy cartesian.

It is sufficient to show that the bisimplicial set map
$\xymatrix@1{\pi^{-1}(\tau)\ar[r]^-{(\zeta_1,\zeta_2)_\ast}
&\pi^{-1}(\sigma)}$ is a pointwise weak equivalence by Proposition
IV.1.7 in \cite {G-J}.

Every bisimplex $\sigma: \triangle^{k,l} \to BC$ is determined by
a string of arrows
$$\xymatrix@1{\sigma: a_0 \ar[r]^-{\alpha_1} &a_1
\ar[r]^-{\alpha_2} &a_2 \cdots \ar[r]^-{\alpha_k} &a_k}$$ of
length $k$ in $C_l$, where $C_l$ is the category in simplicial
degree $l$ in the simplicial category $C$. In horizontal degree
$n$, this bisimplex determines a simplicial set map
$$\bigsqcup_{\gamma: {\bf n} \to {\bf k}}\triangle^l \to BC_n = \underset {(c_0,c_1,\cdots,c_n)}{\bigsqcup}
C(c_0,c_1) \times \cdots \times C(c_{n-1},c_n).$$ On the summand
corresponding to $\gamma: {\bf n} \to {\bf k}$, this map restricts
to the composite
$$\gamma^\ast(\sigma): \triangle^l \to
C(a_{\gamma(0)},a_{\gamma(1)}) \times \cdots \times
C(a_{\gamma(n-1)}, a_{\gamma(n)}) \to BC_n.$$ The simplicial set
$(\underrightarrow{\mbox{holim}}_C X)_n$ in horizontal degree $n$
has the form
$$(\underrightarrow{\mbox{holim}}_C X)_n = \underset{(c_0,c_1,\cdots,c_n)}{\bigsqcup} X(c_0)
\times C(c_0,c_1) \times \cdots \times C(c_{n-1},c_n)$$ It follows
that (in horizontal degree $n$) there is an identification
$$(\pi^{-1}(\sigma))_n = \bigsqcup_{\gamma: {\bf n} \to {\bf k}} X(a_{\gamma(0)}) \times
\triangle^l.$$ It's obvious that any map $(1,\theta):
\triangle^{k,r} \to \triangle^{k,l}$ induces the simplicial set
map
$$\bigsqcup_{\gamma: {\bf n} \to {\bf k}} X(a_{\gamma(0)}) \times
\triangle^r \to \bigsqcup_{\gamma: {\bf n} \to {\bf k}}
X(a_{\gamma(0)}) \times \triangle^l$$ in horizontal degree $n$
which is specified on summands by
$$1 \times \theta: X(a_{\gamma(0)}) \times
\triangle^r \to X(a_{\gamma(0)}) \times \triangle^l$$ such a map
is plainly a weak equivalence, and hence induces a pointwise weak
equivalence
$$\pi^{-1}(\sigma(1\times \theta)) \to \pi^{-1}(\sigma)$$
In particular, any vertex $\triangle^0 \to \triangle^l$ determines
a weak equivalence
$$\bigsqcup_{\gamma: {\bf n} \to {\bf k}} X(a_{\gamma(0)}) \to \bigsqcup_{\gamma: {\bf n} \to {\bf k}}
X(a_{\gamma(0)}) \times \triangle^l$$

Any bisimplicial set map $(\zeta_1, \zeta_2): \triangle^{r,s} \to
\triangle^{k,l}$ and any choice of vertex $v: \triangle^0 \to
\triangle^l$ together induce a commutative diagram of bisimplicial
set maps
$$\xymatrix{
\triangle^{r,0} \ar[r]^-{(1,v)} \ar[d]_-{(\zeta_1, 1)}
&\triangle^{r,s}\ar[d]^-{(\zeta_1, \zeta_2)}\\
\triangle^{k,0} \ar[r]_-{(1, \zeta_2(v))} &\triangle^{k,l}}
$$
It therefore suffices to assume that all diagrams of bisimplices
$$
\xymatrix{ \triangle^{r,0}\ar[dr]^-{\tau}\ar[dd]_\theta\\
&BC\\
\triangle^{k,0}\ar[ur]_-{\sigma}}
$$
induce pointwise weak equivalence
$$\xymatrix@1{\pi^{-1}(\tau) \ar[r]^-{\theta_\ast}&
\pi^{-1}(\sigma)}$$

The simplicial functor $X: C \to {\bf S}$ restricts to an ordinary
functor $X_0: C_0 \to {\bf S}$ via the identification of the
category $C_0$ with a discrete simplicial subobject of the
simplicial category $C$. There is a pullback diagram
$$\xymatrix{
\underrightarrow{\mbox{holim}}_{C_0} X_0 \ar[r] \ar[d] &
\underrightarrow{\mbox{holim}}_C X \ar[d]\\
BC_0 \ar[r] & BC}
$$
and all bisimplices $\triangle^{k,0} \to BC$ factor through the
inclusion $BC_0 \to BC$. Each morphism $a \to b$ of $C_0$ induces
a weak equivalence
$$\xymatrix@1{X_0(a) = X(a)\ar[r]^-\simeq &X(b) = X_0(b)}$$
It therefore follows from the standard argument for ordinary
functors taking values in simplicial sets that all induced maps
$$\xymatrix@1{\pi^{-1}(\tau) \ar[r]^-{\theta_\ast}&
\pi^{-1}(\sigma)}$$ are weak equivalences of simplicial sets.
\end{proof}

\begin{corollary}\label{C:hfib}
Suppose that $G$ is a simplicial groupoid (groupoid enriched in
simplicial sets), and that $X: G \to {\bf S}$ is a simplicial
functor taking values in simplicial sets. Then the map $X(a) \to
F_a$ taking values in the homotopy fibre over $a$ of the diagonal
simplicial set map $\underrightarrow{\mbox{holim}}_G X \to d(BG)$
is a weak equivalence.
\end{corollary}
\begin{proof}
For any arrow $a \to b$ of $G_0$ it has an inverse arrow $b \to
a$, hence the induced simplicial map $X(a) \to X(b)$ has an
inverse map $X(b) \to X(a)$, that means $X(a) \to X(b)$ is a weak
equivalence.
\end{proof}

Suppose that $H$ is an object of the category $s{\bf Gd}$ of
simplicial groupoids and let $f: U \to H$ be a morphism of $s${\bf
Gd}. Take $a \in \mbox{Ob}(H)$ and write $f \downarrow a$ for the
simplicial category given in degree $n$ by the comma category $f_n
\downarrow a$ arising from the functor $f_n: U_n \to H_n$ (Note
that $f \downarrow a$ is a simplicial category in the general
senses, it is not category enriched in simplicial sets since its
simplicial class of objects isn't discrete). Then the functors
$H_n \to cat$ given by $a \mapsto f_n \downarrow a $ assemble to
give a bisimplicial functor $B(f \downarrow~): H \to {\bf S}^2$
with $a \mapsto B(f \downarrow a)$ taking values in bisimplicial
sets. It follows that the assignment $a \mapsto dB(f \downarrow
a)$ defines a simplicial functor taking values in simplicial sets.
The simplicial sets $dB(f \downarrow a)$ therefore become
identified with the homotopy fibres $F_a$ of the diagonal
simplicial set map
$$\xymatrix@1{\underrightarrow{\mbox{holim}}_H dB(f \downarrow ~)\ar[r]^-{d\pi} &dBH}$$
by the Corollary~\ref{C:hfib}. In ``horizontal degree" $n$, this
map can be identified with the projection
$$dB(f \downarrow a_0) \times H(a_0, a_1) \times \cdots \times
H(a_{n-1}, a_n) \to H(a_0, a_1) \times \cdots \times H(a_{n-1},
a_n).$$

The forgetful functors $f_n \downarrow a \to U_n$ also assemble to
define a diagonal weak equivalence
$$\xymatrix@1{\underrightarrow{\mbox{holim}}_H B(f \downarrow~) \ar[r]^-\omega &BU}$$
of trisimplicial sets (here we see
$\underrightarrow{\mbox{holim}}_H B(f \downarrow~)$ as the
trisimplicial set giving rise to the homotopy colimit). This is a
consequence of a standard result of Quillen, applied in each
degree: Suppose that $f: C \to D$ is a functor of small
categories. Then the canonical map
$\underrightarrow{\mbox{holim}}_{D} B(f\downarrow ) \to BC$ is a
weak equivalence. That means the simplicial set map
$$\underrightarrow{\mbox{holim}}_H dB(f \downarrow ~) \to dBU$$
is a weak equivalence.

A \emph{torsor for a presheaf of simplicial groupoids} $G$ is a
simplicial functor $X: G \to {\bf S}Pre(\mathcal{C})$ taking
values in simplicial presheaves such that the associated
simplicial presheaf $\underrightarrow{\mbox{holim}}_G X$ is weakly
equivalent to a point. A morphism $X \to Y$ of $G$-torsors is a
natural transformation of simplicial functors. Insofar as $X$ can
be locally identified with the homotopy fibre of the canonical map
$\underrightarrow{\mbox{holim}}_G X \to dBG$, a map $X \to Y$ of
$G$-torsors restricts to (pointwise) weak equivalences $X|_U \to
Y|_U$ on all sites $ \mathcal{C}\downarrow U$ for which $X(U)$ is
non-empty.

Write {\bf G-Tors} for the category of $G$-torsors, and let
$\pi_0({\bf G-Tors})$ denote its set of path components. There is
a well-defined function
$$\phi: \pi_0({\bf G-Tors}) \to \pi_0({\bf Triv}/dBG)$$
which is induced by associating a $G$-torsor $X$ the element
represented by the map
$$\underrightarrow{\mbox{holim}}_G X \to dBG.$$

There is a function
$$\psi: \pi_0({\bf Triv}/G) \cong \pi_0({\bf Triv}/dBG) \to
\pi_0({\bf G-Tors})$$ which is defined as follows. Let $f: U \to
G$ be an object of ${\bf Triv}/G$ and perform the construction as
above sectionwise to form the diagram
$$\xymatrix{
dBU \ar[d]_-{dBf} & \underrightarrow{\mbox{holim}}_G dB(f
\downarrow
~) \ar[l]_-\simeq \ar[d]^-{f_\ast}\\
dBG & \underrightarrow{\mbox{holim}}_G dB(G \downarrow ~)
\ar[l]^-\simeq_-\alpha \ar[d]^-\beta_-\simeq\\
&dBG}
$$
Then the simplicial $G$-functor $a \mapsto dB(f \downarrow a)$ is
a $G$-torsor since $dBU \to \ast$ is a weak equivalence. This
construction is functorial and defines the function $\psi$.

Note that $\underrightarrow{\mbox{holim}}_G dB(G \downarrow ~)$ is
the simplicial set consisting of strings $(b,a)$ of arrows
$$b_0 \to b_1 \to \cdots \to b_n \to a_0 \to a_1 \to \cdots \to
a_n$$ of length $2n+1$ in $G_n$, and the map $\alpha$ takes this
string to the string $b_0 \to b_1 \to \cdots \to b_n$ while
$\beta$ maps this element to the string $a_0 \to a_1 \to \cdots
\to a_n$, note that $\beta$ is just the projection map $d\pi$. The
$n$-simplices of $\underrightarrow{\mbox{holim}}_G dB(G \downarrow
~)$ can therefore be identified with functors ${\bf n} \ast {\bf
n} \to G_n$ defined on the poset join ${\bf n} \ast {\bf n}$, and
with simplicial structure maps induced by precomposition with maps
$\theta \ast \theta: {\bf m} \ast {\bf m} \to {\bf n} \ast {\bf
n}$. The maps $\alpha$ and $\beta$ are induced by the inclusions
${\bf n} \to {\bf n} \ast {\bf n}$ of the left and right
substrings of length $n$ respectively.

There is a poset map $h_n: {\bf n} \times {\bf 1} \to {\bf n} \ast
{\bf n}$ which is defined by
$$ (i,\epsilon) \mapsto \left\{ \begin{array}{ll}
 i & \mbox{if $\epsilon = 0$, and}\\
n+i & \mbox{if $\epsilon = 1$.}
\end{array}
\right. $$ As a picture, $h_n$ is the diagram
$$\xymatrix{
b_0 \ar[r]\ar[d] &b_1 \ar[r]\ar[d] &\cdots \ar[r] &b_n\ar[d]\\
a_0 \ar[r] &a_1 \ar[r] &\cdots \ar[r] &a_n}
$$

The maps $h_n$ are natural in ordinal number $n$. It follows that
the composites
$$\xymatrix@1{\triangle^n \times \triangle^1 \ar[r]^-{h_n} &B({\bf
n}\ast {\bf n}) \ar[r]^-{(b,a)} &dBG}$$ together define a
simplicial set map $\underrightarrow{\mbox{holim}}_G dB(G
\downarrow ~) \times \triangle^1 \to dBG$ from $\alpha$ to
$\beta$. This construction is natural in all simplicial groupoids,
and so the maps $\alpha$ and $\beta$ are homotopic maps of
simplicial presheaves. Hence the composites $\beta \cdot f_\ast$
and $\alpha \cdot f_\ast$ are homotopic (where $\beta \cdot
f_\ast$ is just the projection map $d\pi:
\underrightarrow{\mbox{holim}}_G dB(f \downarrow ~) \to dBG$). It
follows that the canonical map $\beta \cdot f_\ast$ and the
original map $dBf: dBU \to dBG$ represent the same element of
$\pi_0({\bf Triv}/dBG)$, and so the composite $\phi \cdot \psi$ is
the identity function.

If $X$ is a $G$-torsor, then the canonical map
$\underrightarrow{\mbox{holim}}_G X \to dBG$ is induced by a
simplicial groupoid map $f: E_G X \to G$ (where $E_G X$ is the
simplicial groupoid in which the level $n$ groupoid is the
translation category for the functor $X_n: G_n \to {\bf Set}$ in
each degree $n$). The comma category $f_n \downarrow a$ has
objects $((b,x), b \to a)$ where $a, b \in G_n, x \in X_n(b)$ and
morphisms $\alpha: ((b,x), b \to a) \to ((c,y), c \to a)$ where
$\alpha: b \to c$ is a morphism of $G_n$ such that $X_n(\alpha)(x)
= y$ and the composite with $c \to a$ is the map $b \to a$. The
$n$-simplices in $dB(f \downarrow a)$ is:
$$((b_0, x_0), b_0 \to a) \to ((b_1, x_1), b_1 \to a) \to \cdots
\to ((b_n, x_n), b_n \to a)$$ the $n$-simplices in $B(f_n
\downarrow a)_n$. There is a $G$-natural function $dB(f \downarrow
a)_n \to X_n(a)$ sending the $n$-simplex to $X_n(\beta)(x_0)$
where $\beta: b_0 \to a$. For every $n$-simplex $x$ in $X_n(a)$
its preimages are connected, hence these functions induce a map
$dB(f \downarrow a) \to X(a)$ for all $a$ which is a weak
equivalence, and hence determines a map of simplicial functors
$$ dB(f \downarrow ~) \to X$$
$X$ is a $G$-torsor, then $\underrightarrow{\mbox{holim}}_G X$ is
weakly equivalent to a point, so is
$\underrightarrow{\mbox{holim}}_G dB(f \downarrow~)$, that means
$dB(f \downarrow~)$ is a $G$-torsor as well, hence the above map
is a map of $G$-torsors. It follows that the composite $\psi \cdot
\phi$ is the identity function. We have therefore proved the
following

\begin{theorem}\label{T:maintor}
The natural function
$$\phi: \pi_0({\bf G-Tors}) \to \pi_0({\bf Triv}/dBG) \cong [\ast, \overline{W}G] $$
is a bijection for each presheaf of simplicial groupoids $G$ on
any Grothendieck site $\mathcal{C}$.
\end{theorem}
\begin{proof}
The displayed isomorphism follows from Lemma \ref{L:Diagram} and
Lemma \ref{L:Triv2}. The proof that $\phi$ is a bijection is
displayed above.
\end{proof}

\begin{remark}
The definition of $G$-torsors and the bijection in Theorem
\ref{T:main} are available for the sheaf of simplicial groupoids
$G$. If G is a presheaf of simplicial groupoids and $L^{2}G$ is
its associated sheaf, then $\overline{W}L^{2}G$ is the simplicial
sheaf associated to $\overline{W}G$, and the map $\overline{W}G
\to \overline{W}L^{2}G$ is a weak equivalence, and all we're
interested in is the invariant $[\ast,\overline{W}G] =
[\ast,\overline{W}L^{2}G]$.

Same arguments are valid in Section 3.3 and 3.4.
\end{remark}

\begin{remark}
Joyal-Tierney \cite {J-T6} obtain similar result for the sheaves
of simplicial groupoids $G$, but their definition of $G$-torsors
is different from ours, our definition is much more flexible.
\end{remark}

Theorem~\ref{T:2gtor} is a special case of
Theorem~\ref{T:maintor}.

\begin{example}
When a Grothendieck site $\mathcal{C}$ is the trivial category
$\ast$, a presheaf of simplicial groupoids $G$ over $\mathcal{C}$
is just an ordinary simplicial groupoids. Thus $\pi_0({\bf
G-Tors}) \cong [\ast, \overline{W}G] \cong [\ast, G] \cong \pi_0
G$, so the set of path components of $G$-torsors is bijective to
the set of path components of the simplicial groupoid $G$.
\end{example}

\bigskip
\noindent{\bf Acknowledgement}

\noindent I would like to thank my supervisor Dr. J.F. Jardine for
suggesting the topic, and for his help and continuing
encouragement.


\begin{thebibliography}{99}

\bibitem{B} L. Breen, \emph{On the classification of 2-gerbes
and 2-stacks}, Ast$\acute{\mbox{e}}$risque 225, 1994.

\bibitem{BK} K.S. Brown, \emph{Abstract homotopy theory and
generalized sheaf cohomology}, Trans. AMS {\bf 186}, 419-458,
1973.

\bibitem{B-G} K.S. Brown and S.M. Gersten, \emph{Algebraic
K-theory as generalized sheaf cohomology}, Lecture Notes in Math.
341, Springer-Berlin, 266-292, 1973.

\bibitem{Bry} J. Brylinski, \emph{Loop spaces, Characteristic
classes and geometric quantization}, Birkh$\ddot{a}$user, PM107,
1993.

\bibitem{Cr} S.E. Crans, \emph{Quillen closed model structures
for sheaves}, J. Pure Appl. Algebra 101:35-57, 1995.

\bibitem{D-K} W.G. Dwyer and D.M. Kan, \emph{Homotopy theory and
simplicial groupoids}, Indag. Math. $\textbf{46}$, 379-385, 1984.

\bibitem{G} J.Giraud, \emph{Cohomologie non ab$\acute{e}$lienne}, Grundlehren
$\#$179, Springer Verlag, 1971.

\bibitem{G-J} P.G. Goerss and J.F. Jardine, \emph{Simplicial Homotopy Theory},
Birkh$\ddot{\mbox{a}}$user, PM174, 1999.

\bibitem{J1}  J. F. Jardine, \emph{Simplicial presheaves},
J. Pure and Appl. Algebra {\bf 47}, 35-87, 1987.

\bibitem{J3}  J. F. Jardine, \emph{Universal Hasse-Witt classes},
Contemporary Math. {\bf 83}, 83-100, 1989.

\bibitem{J4}  J. F. Jardine, \emph{The homotopical foundations of
algebraic K-theory},  Contemporary Math. {\bf 83}, 57-82, 1989.

\bibitem{J5}  J. F. Jardine, \emph{Stacks and the homotopy theory of simplicial sheaves},
 Homology, Homotopy and Applications 3(2), 361-384, 2001.

\bibitem{J6}  J. F. Jardine, \emph{Notes}, 2003.

\bibitem{J7}  J. F. Jardine, \emph{Simplicial objects in a
Grothendieck topos}, Contemporary Math. 55 (1), 193-239, 1986.

\bibitem{Jo} A. Joyal, \emph{letter to A.Grothendieck}, (1984).

\bibitem{J-T}  A. Joyal and M. Tierney, \emph{Strong stacks and classifying
space}, Category Theory (Como, 1990), Springer Lecture Notes in
Math. {\bf 1488} (1991), pp 213-236.

\bibitem{J-T5}  A. Joyal and M. Tierney, \emph{On the homotopy
theory of sheaves of simplicial groupoids}, Math. Proc. Camb.
Phil. Soc. {\bf 120}, pp 263-290, 1996.

\bibitem{J-T6}  A. Joyal and M. Tierney, \emph{Classifying spaces for sheaves
of simplicial groupoids,} Journal of pure and applied algebra {\bf
89}, pp 135-161, 1993.

\bibitem{Luo} Z. Luo, \emph{Closed model categories for presheaves of simplicial
groupoids and presheaves of 2-groupoids}, Preprint ArXiv:
math.AT/0301045.

\bibitem{Luo1} Z. Luo, \emph{Presheaves of simplicial
groupoids}, Ph.D thesis, University of Western Ontario, 2004.

\bibitem{Mac} S. Mac Lane, \emph{Categories for the working
mathematician}, second edition, Graduate texts in Mathematics 5,
Springer, 1998.

\bibitem{M-M} S. Mac Lane and I. Moerdijk, \emph{Sheaves in geometry and logic}, Springer-Verlag, 1992.

\bibitem{M-S}  I. Moerdijk and J. Svensson, \emph{Algebraic classification of
equivariant homotopy 2-types, I}, J. Pure Appl. Algebra
89:187-216, 1993.

\bibitem{M-V} F. Morel and V. Voevodsky, \emph{$ A\sp 1$-homotopy theory of
schemes}, Inst. Hautes Études Sci. Publ. Math. No. 90, 45--143,
1999.


\bibitem{Q1} D. Quillen, \emph{Homotopical Algebra}, Lecture Notes in
Math., Vol.43, Springer, Berlin-Heidelberg-New York, 1967.

\bibitem{Q2} D. Quillen, \emph{Rational homotopy theory}, Ann. of
Math. $\bf 90$, 205-295, 1969.


\end{thebibliography}
\end{document}